\documentclass[12pt,a4]{article}
\usepackage{geometry}
\usepackage{multirow}
\usepackage{graphicx}
\usepackage{graphics} % for pdf, bitmapped graphics files
\usepackage{epsfig}
\usepackage{amsmath, amsthm, amssymb, amsfonts, enumerate}
\usepackage{blindtext}
\usepackage{float}
\usepackage{colortbl}
\usepackage[english]{babel}
\usepackage{authblk}
\usepackage[colorlinks=true, pdfstartview=FitV,breaklinks=true, 
linktocpage = true, citecolor = magenta]{hyperref}
\usepackage{algorithm}

\usepackage[noend]{algpseudocode}
\usepackage{times} % assumes new font selection scheme installed
\usepackage{url}
\usepackage{booktabs}
\usepackage{array}

\usepackage{sidecap}
\sidecaptionvpos{figure}{t}

\usepackage[toc,page]{appendix}

\usepackage{subcaption}

\usepackage[usenames,dvipsnames]{pstricks}
\usepackage{pst-grad} 
\usepackage{threeparttable}

\newtheorem{lemma}{Lemma}
\newtheorem{proposition}{Proposition}
\newtheorem{theorem}{Theorem}
\newtheorem{corollary}{Corollary}

\theoremstyle{definition}
\newtheorem{definition}{Definition}
\theoremstyle{definition}
\newtheorem{remark}{Remark}
\theoremstyle{definition}

\algrenewcommand\algorithmicrequire{\textbf{Input:}}
\algrenewcommand\algorithmicensure{\textbf{Output:}}

\usepackage{lineno}

\begin{document}

	\title{Finding Approximately Convex Ropes \\\medskip in the Plane}

	% for line numbering
	%\linenumbers

	\author[1,4]{Le Hong Trang}
	\author[3]{Nguyen Thi Le}
	\author[2,4]{Phan Thanh An \thanks{Corresponding author. E-mail: \texttt{thanhan@hcmut.edu.vn}}}
	%%\thanks{\it Corresponding author: thanhan@math.ac.vn},
	\affil[1]{\small Faculty of Computer Science and Engineering, Ho Chi Minh City University of Technology (HCMUT), 268 Ly Thuong Kiet Street, District 10, Ho Chi Minh City, Vietnam}
	\affil[2]{\small Institute of Mathematical and Computational Sciences (IMACS) and Faculty of Applied Science,
		Ho Chi Minh City University of Technology (HCMUT), 268 Ly Thuong Kiet Street, District 10, Ho Chi Minh City, Vietnam}
	\affil[3]{\small Institute of Mathematics,  Vietnam Academy of Science and Technology (VAST), 18 Hoang Quoc Viet Road, %Cau Giay District, 
		Hanoi, Vietnam}
	\affil[4]{\small  Vietnam National University Ho Chi Minh City, Linh Trung Ward, Thu Duc City, Vietnam}
	\date{}
	
	\maketitle
	\newenvironment{mylisting}
	{\begin{list}{}{\setlength{\leftmargin}{1em}}\item\normalsize\bfseries}
		{\end{list}}
	\newenvironment{mytinylisting}
	{\begin{list}{}{\setlength{\leftmargin}{1em}}\item\scriptsize\bfseries}
		{\end{list}}
	
	% 	\author[1,2]{\fnm{Le Hong Trang}\sur{}}%{\email{lhtrang@hcmut.edu.vn}}
	% 	\author[3,4]{\fnm{Nguyen Thi Le}\sur{}}%{\email{ngle37@gmail.com}}
	% 	\author*[1,2]{\fnm{Phan Thanh An}\sur{}}%{\email{thanhan@hcmut.edu.vn}}

	%
	
	\begin{abstract}
		The convex rope problem
		is to find a counterclockwise 
		or clockwise convex rope starting at the vertex $a$ and ending at the vertex $b$ of a simple
		polygon $\mathcal{P}$, where $a$ is a vertex of the convex hull of $\mathcal{P}$ and  $b$ is 
		visible from infinity. The  convex rope mentioned 
		is   the shortest path joining $a$ and $b$ that does not enter the interior of $\mathcal{P}$.
		In this paper, the problem is  reconstructed  as the one of finding such shortest path in a simple polygon and solved by the method of multiple shooting.
		We then show that  if the collinear condition of the method holds at all shooting points, then these shooting points form the shortest  path.  Otherwise, 
		the sequence of paths obtained  by the update of the method converges to the  shortest path. 
		The algorithm is implemented in C++ for numerical experiments.

		
		\medskip
		{\bf MSC2010:} MSC 52A30; MSC 52B55; MSC 68Q25; 90C26; MSC 65D18; 68R01.
		
	\end{abstract}
	
	%	\keywords{Convex hulls, approximate algorithm, convex rope, shortest path, non-convex optimization}
	
	%	\pacs[MSC2010]{ MSC 52A30; MSC 52B55; MSC 68Q25; 90C26; MSC 65D18; 68R01.}
	%
	
	\maketitle
	
	\section{Introduction}
	\label{sec:intro}
	%An autonomous grasp
	%planning for robotic hands, the convex rope problem, posted by Peshkin and Sanderson in~\cite{peshkin} plays a significant role in  planning reachable grasps of a  polygonal object by a simple
	%robot hand. 
	Solving  automatic grasp planning problems is to  evaluate must-touch regions, what forces should be applied to an object, and how those forces can be used by robotic hands~\cite{perez, SongAtAl2018, Xue2008}.
	A geometric construction, namely the convex rope posted by Peshkin and Sanderson in 1986~\cite{peshkin} gives valuable information for  grasping automatically with robot hands. The problem of finding convex ropes plays a significant role in grasp planning  a  polygonal object with a simple robot hand. 
	Besides the weight information of the object relevant to grasping and their geometric information  are also necessary for this task. To simplify and decrease the number of possible hand configurations, objects  are modeled as a set of shape primitives, such as polygons in 2D or polyhedrons, spheres, cylinders, and cones in 3D. 
	
\begin{figure}
		\centering
		\begin{subfigure}[b]{0.3\textwidth}
			\centering
			\includegraphics[scale=0.15]{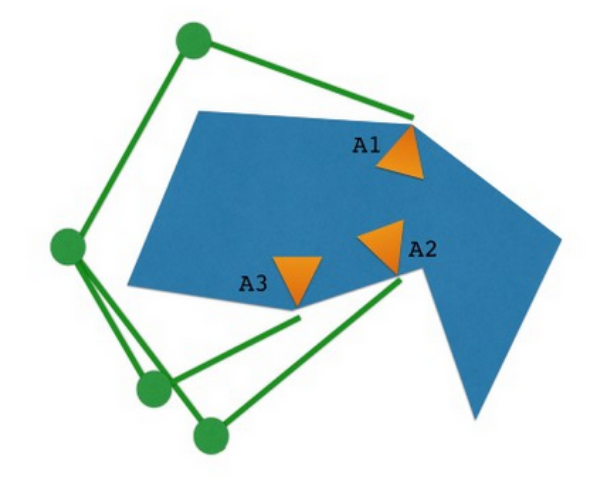}
			\caption*{(i)}
		\end{subfigure}
		\begin{subfigure}[b]{0.3\textwidth}
			\centering
			\includegraphics[scale=0.2]{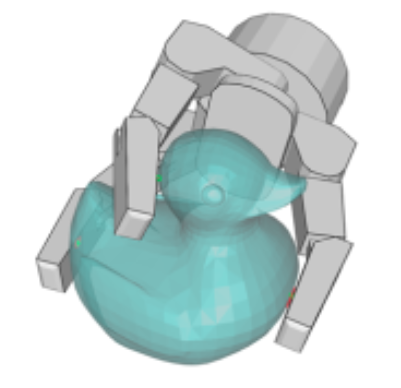}
			\caption*{(ii)}
		\end{subfigure}
		\begin{subfigure}[b]{0.3\textwidth}
			\centering
			\includegraphics[scale=0.135]{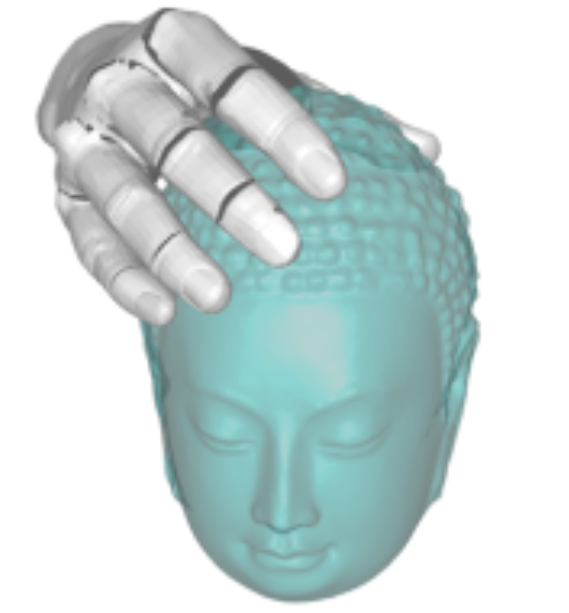}
			\caption*{(iii)}	
		\end{subfigure}
		\caption[]{Illustration of  grasping with a simple robot hand; 	(i): grasping  a polygonal object in 2D, images from \textit{https://www.shadowrobot.com}; (ii) and (iii): grasping  a polyhedral object in 3D,  images from~\cite{SongAtAl2018}.}
		%
		%	\caption[]%
		%	{Real caption\footnote{blah}}
		\label{fig:applications}
	\end{figure}

\begin{figure}
		\centering
		\begin{subfigure}[b]{0.4\textwidth}
			\centering
			\includegraphics[scale=0.5]{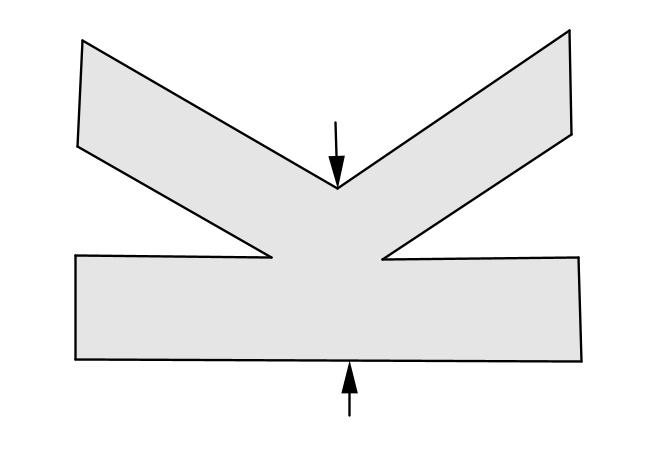}
			%	\caption*{(a)}
		\end{subfigure}
		\begin{subfigure}[b]{0.4\textwidth}
			\centering
			\includegraphics[scale=0.4]{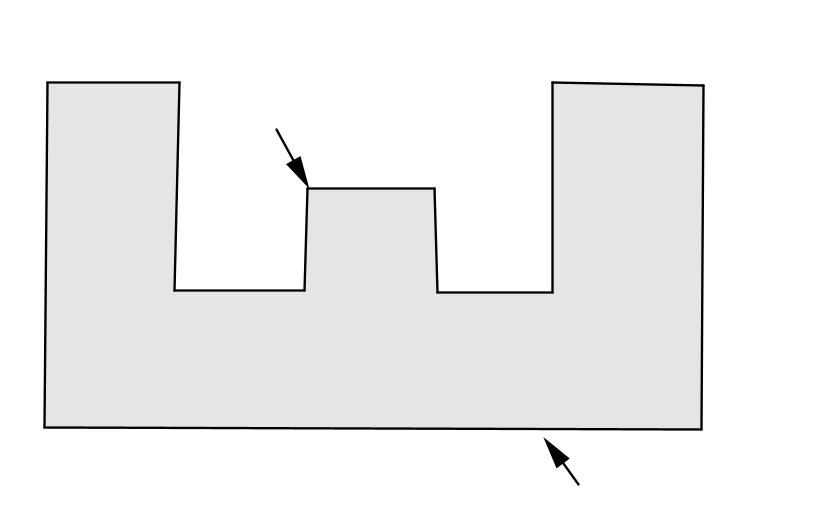}
			%	\caption*{(b)}
		\end{subfigure}
		\caption{Contact points would not necessarily be in edges of the convex hull of the polygon.}
		\label{fig:letter}
	\end{figure}

	In 2D, determining  convex ropes helps to find contact points  which are  points on a polygonal object corresponding to the robot's  finger joint position (see Fig.~\ref{fig:applications}(i)). For an object in 3D, the automatic grasp planning  concentrates on finding must-touch regions on the object to get information about possible hand configurations (see Figs.~\ref{fig:applications}(ii) and (iii)).
	Contact points usually belong to the trajectory of the convex rope joining two vertices of a simple polygon. These two vertices are not required to be located on the edges of the convex hull of the polygon (see Fig.~\ref{fig:letter}). 
	It turns out to consider that the endpoints of the convex rope can be visible from infinity (see in Sect.~\ref{sect:problem-approach} for the definition).
	%The problem can be also stated as a non-convex optimization problem as follows
	In the non-convex optimization context, the problem can be stated as follows
	\begin{align*}
		%	\begin{split}
		\min_{(k, x_0, x_1, \ldots, x_{k+1})} & \sum_{i = 0}^{k} \|x_i-x_{i+1}\| \\
		\text{subject to}\, & x_i \in \mathcal{P}, \text{for}\, i = 1, 2, \ldots, k, \, x_0 = a, x_{k+1} = b\\
		&]x_i,x_{i+1}[\cap\hbox{int}(\mathcal{P})=\emptyset, \text{ where } ]x,y[ := [x,y] \setminus \{ x, y \}.
		%	\end{split}
	\end{align*}
	Unfortunately, no optimization methods have been found for solving the problem.
	In this paper, we consider geometrical methods for solving the problem to get 
	the global solution. 
	
	%It stands to reason that we can consider the case that endpoints of the convex rope are visible from infinity (its  detailed definition will appear in Sect.~\ref{sect:problem-approach}).
	%
	%Algorithms solving the convex rope problem apply directly only to  polygons. For 3D objects, when their 2D cross-section contains all the information about the object,  the convex rope construction and the algorithm can be applied to the case.
	%
	%The grasping problem is to find surfaces on the object, and a configuration of a robot’s fingers, satisfying three conditions:
	%\begin{itemize}
	%	\item The fingers must be in contact with the object.
	%	\item The configuration must be reachable, i.e., there must 	be a collision-free path for the robot to get to the 	grasping configuration.
	%	\item The object must be stable once grasped, i.e., it must 	not slip during subsequent motion.
	%\end{itemize}

	%\footnotetext[2]{source of the images from~\cite{SongAtAl2018}}
	
	% Motivation of finding convex ropes and its application were presented carefully in [???].
	%\textit{Related Work}: 
	%
	%To date, to solve the convex rope problem, there have been two main viewpoints: convex hull approach and shortest path approach.
	The non-convex optimization problem can be addressed via solving alternatively one of two fundamental problems which are convex hull and shortest path problems.
In 1986 an algorithm based on evaluating cumulative angles  was introduced~\cite{peshkin}, in which these angles are calculated for each side of the polygon.  Another algorithm  proposed by An in 2010~\cite{ptan-convex-rope} was derived from Melkman convex hull algorithm.
Using triangulation, a linear time algorithm in 1987 was presented in~\cite{guibas-funnel}
to compute convex ropes as boundaries of the convex hull of a polyline.  In 1990, Heffernan and Mitchell~\cite{heffernan} also addressed a linear time algorithm by sorting to a partial triangulation of a polygon. 
All  algorithms mentioned are exact. In 2006, Li and Klette~\cite{Li2006} introduced an approximate algorithm 
%(say Rubber band algorithm) 
that starts with a special trapezoidal segmentation of a polygon  and  their segmentation is simpler than the triangulation procedure. 
%This algorithm does not depend on a  rather complicated, but linear-time triangulation. 
%
%In the paper, we also present an approximate algorithm based on the multiple shooting approach   for determining  convex ropes.

In this paper, we present an approximate algorithm using the geometric multiple shooting approach, which was proposed for solving the geometric shortest path problems in 2D and 3D~\cite{AnHaiHoai2013, AnTrang2018, HoaiAnHai2017}.
In particular, the convex rope problem will be formulated into the form of finding shortest paths in a simple polygon since the closed relationship between convex sets and shortest paths. A set is convex if the line segment joining any two points of the set lies completely inside the set.
Substituting a line segment joining two points by the shortest path joining these two points expands the notion of classical convexity to
geodesic convexity~\cite{HaiAn2011,Toussaint1989}.
	We then use successfully the method of multiple shooting (MMS for short) to find approximately these shortest paths. 
	The three main factors of the method are given in detail.
	As an advantage of MMS, instead of solving the problem for the whole formed polygon, we deal with a set of smaller partitioned subpolygons. This is useful in the implementation of devices with limited memory. 
	%
	%In particular, the counterclockwise (clockwise, respectively) convex ropes can be obtained by  approximate shortest paths inside  other polygons which are determined using the idea of the method of multiple shooting in geometry
	%~\cite{bock-dmsm}. 
	%
	%
	%Let us recall some previous works related to the multiple shooting method. For solving   ordinary differential equation  boundary value problems, a numerical method called the method of multiple shooting (MMS for short) is introduced.  Based on the idea, MMS  in computational geometry  has been posed in 2013 (see~\cite{AnHaiHoai2013}) by An, Hai, and Hoai   for finding approximately  shortest paths joining two given points in a  simple polygon. They expand the multiple shooting approach   for solving approximately  shortest path problem on  convex polytopes in 3D (see~\cite{AnTrang2018, HoaiAnHai2017}). 
	%
	%
	%%%%%%%%%%%%%%  Must not detail such as  2021-5-11%%%%%%%%%%%%%%%%%%%%%
	
	In~\cite{AnHaiHoai2013}, the multiple shooting approach of An et al.  is applied to  simple polygons  with a  hypothesis of the so-called ``oracle'', not for general simple polygons. The ``oracle'' condition requests that between two consecutive cutting segments there exists at least one point in which every shortest path
	joining two corresponding shooting points passes through it, where  notions of cutting segments and    shooting points are introduced in Sect.~\ref{sect:MMS}.                                                                                                                                                                                                                                  
	Besides,  the relationship between the path satisfying the stop condition  and the shortest path; the convergence of paths obtained by  An et al.'s algorithm is not stated. 
	In the sequence, there are some open questions related to the algorithm in~\cite{AnHaiHoai2013} as follows
	\begin{itemize}
		\item Can use MMS for  simple polygons that are not ''oracle''?
		\item If   $Z^*$ is a path obtained by the algorithm at an iteration which satisfies the stop condition, whether $Z^*$ is the shortest path or not?
		\item Does the sequence of paths $\{ Z^j\}_{j=0}^{+\infty}$  obtained by the algorithm  converges to the  shortest path?
	\end{itemize}
	%%%%%%%%%%%%%%  Must not detail such as%%%%%%%%%%%%%%%%%%%%%
	%
	%\textit{Our contributions}: 
	%In the paper, we use successfully MMS to deal with approximately the convex rope problem. 
	%The three main factors of the method are given in detail. 
	%
	%	Hai and An in~\cite{HaiAn2011} showed the relationship between  types of convergence of  sequences of paths in a simple polygon,  namely that with respect to the Hausdorff distance and the length.
	
	In simple polygons,  Hai and An in~\cite{HaiAn2011} showed the relationship between the convergence of  sequences of paths  with respect to the Hausdorff distance and the one  with respect to  the length.
	In this paper, we will apply this result to get the convergence  to the shortest path of the sequence of paths obtained by the proposed algorithm. 		 
	Particularly, we adapt MMS to the convex rope problem with slight modifications on determining  the stop condition and updating new paths to answer these open questions.%
	We construct a simple polygon in which our corresponding algorithm  can  apply even though the polygon does not satisfy the ``oracle'' condition (Sect.~\ref{sect:creating-polygon}). 
	%	It is known that using the same approach in 3D as in~\cite{AnTrang2018, HoaiAnHai2017} only gets local solutions while 
	We can also show that, at an iteration step, if a path satisfies the stop condition at all shooting points, then the path and the shortest path are identical. Otherwise, Theorems~\ref{theo:global solution} and \ref{theo:update_path} state that the sequence of paths obtained by the algorithm  converges to the shortest path, then the global solution is obtained. 
	
	%
	%	Unlike previous works using MMS, we here show that the sequence of paths obtained by our algorithm converges to the  shortest path by the Hausdorff distance (Theorem~\ref{theo:global solution}), then the global solution is obtained. 
	%	In addition, we verify that  at an iteration step, if a path satisfies the stop condition at all shooting points, then the path and the shortest path are identical (Proposition~\ref{prop:collinear-condition}). Since the length function is lower semi-continuous but not continuous, the convergence of the sequence of paths obtained by the algorithm does not guarantee the convergence of the sequence of lengths of these paths. Therefore Theorem~\ref{theo:update_path} states that the sequence of lengths of paths obtained by the algorithm are strictly decreasing and convergent if the stop condition does not hold for at least one shooting point.
	
	\medskip
	The rest of the paper is organized as follows. Sect.~\ref{sect:problem-approach} introduces the convex rope problem and three main factors of MMS.
	%  recalls geometric preliminaries.  
	Sect.~\ref{sect:MMS} presents the proposed algorithm
	% Here the convex rope problem is solved by referring to the problem of finding the shortest path joining two points in a simple polygon. 
	%	The construction of the polygon is shown in Sect.~\ref{sect:creating-polygon} 
	in which MMS with three factors (f1)-(f3) can be applied.
	% in Sect.~\ref{sect:partition},~\ref{sect:collinear_condition} and~\ref{sect:update}. 
	The correctness of the proposed algorithm 
	and answers for the above open questions 
	%improved results
	%
	is stated by results in Sect.~\ref{sect:correctness}.
	%
	%	We  show that if the colinear condition  holds at all shooting points, then these shooting points form the shortest  path (Proposition~\ref{prop:collinear-condition}). Otherwise,  	the sequence of paths obtained by the algorithm converges, by means of the Hausdorff distances, to the  shortest path (Theorem~\ref{theo:global solution}). 
	%the sequence of length of paths  found by the Algorithm~\ref{alg:msm-convex_rope} is strictly decreasing and convergent (Theorem~\ref{theo:update_path} and Corollary~\ref{coro:convergent}).
	The algorithm is implemented in C++ and numerical results are given
	and visualized in Sect.~\ref{sect:implementation} to describe how our method works.
	Some advantages of  MMS are established in Sect.~\ref{sect:conclusion}. Geometrical properties and
	their proofs  are arranged in Appendix. % (Sect.~\ref{sect:appendix}).
	
	\section{Preliminaries}
	\label{sect:problem-approach}
	
\begin{figure}[htp]
		\centering
		\begin{subfigure}[b]{0.5\textwidth}
			\centering
			\includegraphics[scale=0.32]{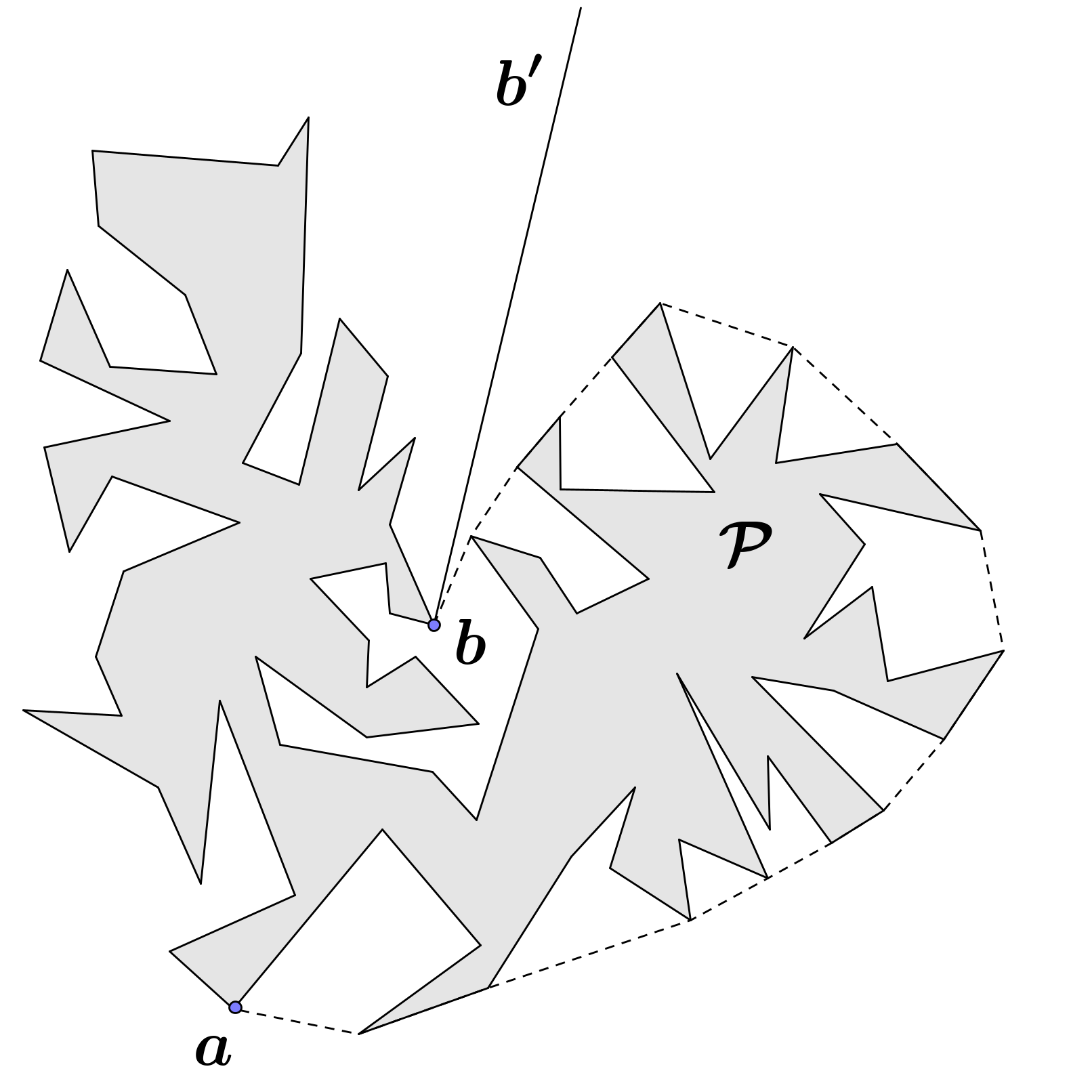}
			\caption*{(i)}
		\end{subfigure}%
		\begin{subfigure}[b]{0.5\textwidth}
			\centering
			\includegraphics[scale=0.31]{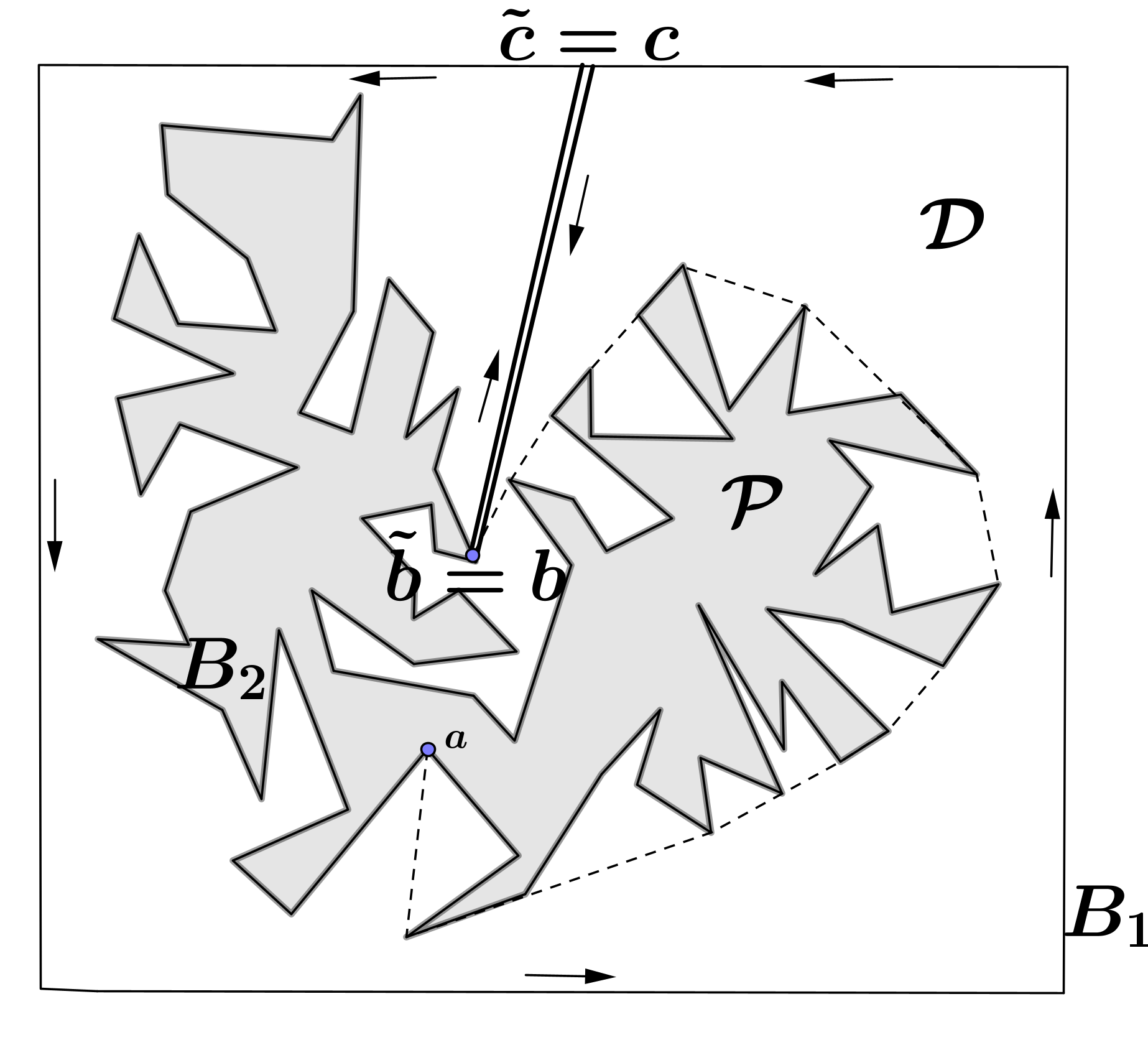}
			\caption*{(ii)}
		\end{subfigure}
		\caption{The counterclockwise convex rope starting at $a$ and ending at $b$, the case (i): $a$ is a vertex of the convex hull of $\mathcal{P}$ and the case (ii): $a$ is visible from infinity.}
		\label{fig:construct-polygon}
	\end{figure}

	A vertex of a simple polygon $\mathcal{P}$ is \textit{visible
		from infinity} if there exists a ray beginning at the vertex which does not
	intersect $\mathcal{P}$ anywhere else. In Fig.~\ref{fig:construct-polygon}(i), $b$ is
	visible from infinity with ray $bb'$. Clearly, a vertex of the convex hull of $\mathcal{P}$ is also visible from infinity.
	We say a path starting at the vertex $a$ and ending at the vertex $b$ of  $\mathcal{P}$ is \textit{counterclockwise} (\textit{clockwise}, resp.) if
	the interior of $\mathcal{P}$ lies to the left (right, resp.) as we step
	along the path, starting at $a$. 
	\begin{definition}
		The counterclockwise (clockwise, resp.)
		shortest path starting at $a$ and ending at  $b$ of
		$\mathcal{P}$ that does not enter the interior of $\mathcal{P}$, where $a$ is a vertex of the
		convex hull of $\mathcal{P}$ and $b$ is visible from infinity, is called a
		\textit{counterclockwise (clockwise, resp.) convex rope} starting at $a$ and ending at  $b$. 
	\end{definition}

	The convex rope problem  is necessary to find a counterclockwise (clockwise, resp.)
	convex rope between $a$ and  $b$.
	For simplicity, we consider the
	counterclockwise convex rope case only. The clockwise convex rope case is
	considered similarly.  
	%Furthermore, we  	assume that the location  of $a, b$, and the ray $bb'$ is given in Fig.~\ref{fig:construct-polygon}(i) in which the counterclockwise convex rope between $a$ and $b$ is illustrated as dashed curve.
	%
	The counterclockwise convex rope between $a$ and $b$ is illustrated as a dashed curve  in Fig.~\ref{fig:construct-polygon}(i).
	When both  $a$ and $b$ are visible infinity, we obtain another version of the convex rope problem which can be also solved by the same way as the case that $a$ is a vertex of the convex hull of $\mathcal{P}$  by our algorithm. For details see Sect.~\ref{sect:creating-polygon}.
	%	\begin{figure}[htp]
	%		\centering
	%		\includegraphics[width=0.33\linewidth]{overview_0}
	%		\caption{The counterclockwise convex rope starting at $a$ and ending at $b$.}
	%		\label{fig:overview0}
	%	\end{figure}

	% $[p, q] := \{(1-\lambda)p + \lambda q : 0 \le \lambda \le 1\}, (p, q) := \{(1-\lambda)p + \lambda q : 0 < \lambda < 1\} $.
	
	% For all $x \in ]p,q[$, $x$ is called an \textit{interior point} of  $[p,q]$.
	%Given three vertices $p = (p_x, p_y), q = (q_x, q_y)$, and $r = (r_x, r_y)$, let $s =(r_x - p_x)(p_y - q_y) - (q_x - p_x)(p_y - r_y)$. We say that $r$ is \textit{left} of (respectively is on, is \textit{right} of) the ray $pq$ (an oriented line from $\gamma$ to $q$) if $s$ is positive (respectively zero, negative).
	%
	
	We now recall some basic concepts and properties. For any points $x, y$ in the plane, we denote 
	$[x,y] := \{ (1-\lambda)x + \lambda y: 0 \leq \lambda \leq 1 \}$, 
	$]x,y] := [x,y] \setminus \{ x \}$, 
	$]x,y[ := [x,y] \setminus \{ x, y \}.$ 
	%For any points $x$ and $y$ in a simple polygon, there is a \textit{shortest path} joining $x$ and $y$ in the polygon, denoted by $\text{SP}(x,y)$. 
	Let $x$ and $y$ be two points in a simple polygon $\mathcal{A}$. The \textit{shortest path} joining $x$ and $y$ in $\mathcal{A}$, denoted by $\text{SP}(x,y)$, is the path of minimal length joining $x$ and $y$ in $\mathcal{A}$. As shown in \cite{Chazelle1982,Lee1984}, $\text{SP}(x,y)$ is a polyline whose vertices except for $x$ and $y$ are reflex vertices of $\mathcal{A}$, i.e., the internal angles at these vertices are at least $\pi$. In addition, these vertices are also reflex vertices of $\text{SP}(x,y)$. Furthermore, let $y' \in \mathcal{A}$, for any two line segments $e$ of $\text{SP}(x,y)$ and $g$ of $\text{SP}(x,y')$, we have either (i) $e=g$; or (ii)  $e$ and $g$ are disjoint or share at most one endpoint.
	%
	% 	\begin{remark}[\cite{Chazelle1982,Lee1984}] We have
	% 		\label{rem:shortest_path}
	% 		\begin{itemize}
	% 			\item[(a)] 	$\text{SP}(x,y)$ is a polyline whose vertices except for $x$ and $y$ are reflex vertices of $\mathcal{A}$, i.e., the upper angles at these vertices are at least $\pi$. In addition, these vertices are also reflex vertices of $\text{SP}(x,y)$.
	% 			\item[(b)] Let $y' \in \mathcal{A}$. For any two line segments $e$ of $\text{SP}(x,y)$ and $g$ of $\text{SP}(x,y')$, we have either (i) $e=g$; or (ii)  $e$ and $g$ are disjoint or share at most one endpoint.
	% 		\end{itemize}
	% 	\end{remark}
	Fundamental concepts such as polylines, polygons, convex (reflex) vertices, shortest paths and their properties which will be used in this paper can be seen in~\cite{An2017,HaiAnHuyen2019} and Appendix.
	%%%

	Let $x \in \mathbb{R}^2$ and $A$ be a nonempty subset of $\mathbb{R}^2$, then the distance from $x$ to $A$, denoted by $d (x, A)$, is defined as 	$d (x, A) = \inf_{ y \in A} \Vert x-y \Vert,$
	where $\Vert.\Vert$ is the  Euclidean  norm in $\mathbb{R}^2$. 
	Let $A$ and $B$ be two nonempty subsets of $\mathbb{R}^2$. The \textit{Hausdorff distance} between $A$ and $B$, denoted by $d_{\mathcal{H}}(A,B)$ is defined as $d_{\mathcal{H}}(A,B)=\max \{ \sup_{x \in A} d(x,B), \sup_{y \in B} d(y,A) \}$.
	% 	Equivalently, we have
	% 	\begin{align}
	% 		\label{def:Hausdorff-distance}
	% 		d_{\mathcal{H}}(A,B)=\inf \{ \epsilon \ge 0 \text{ such that } B \subset N(A, \epsilon) \text{ and } A \subset N(B, \epsilon)\},
	% 	\end{align}
	% 	where $N(S, \epsilon) = \{x \in \mathbb{R}^2 \vert \, d (x, S) \le \epsilon\},$ for $\epsilon >0$ and $S \subset \mathbb{R}^2.$
	
	\medskip
	Regarding MMS for finding the shortest path joining two points $p$ and $q$ in a domain $\mathcal{D}$ \cite{AnHaiHoai2013,AnTrang2018, HoaiAnHai2017}, three following factors are included:
	\begin{itemize}
		\item[(f1)] Partition of the  domain $\mathcal{D}$ which is a polygon in 2D (polytope in 3D, resp.) into subpolygons (subpolytopes, resp.) is created by cutting segments (cutting slices, resp.).  In each  cutting segment  (the boundary of cutting slice, resp.) we take a point that is called a shooting point;
		\item[(f2)] Construct a path in $\mathcal{D}$ joining  $p$ and $q$, formed by the ordered set of shooting points. 
		% The path is the concatenation of the shortest paths joining two consecutive shooting points  in corresponding subdomains.
		A stop condition called collinear condition (straightness condition, resp.) is established at shooting points;
		\item[(f3)]  The algorithm enforces (f2) at all shooting points to check the stop condition. Otherwise, an update of shooting
		points improves the paths joining  $p$ and $q$.% that are formed by the set of new shooting points. 
	\end{itemize}
	These factors which are established for the convex rope problem  are described in the next section.
	\section{MMS-based Algorithm for Finding  Convex Ropes}
	\label{sect:MMS}

	%To find the convex rope starting at $a$ and ending at $b$, there are two  steps in the propose algorithm as follows:
	The proposed algorithm includes two following phases. We first construct a simple polygon such that the convex rope problem is referred to as the problem of finding the shortest path joining two points in the polygon. In the second phase, three factors (f1) - (f3) will be applied for the shortest path problem. 
	%\begin{itemize}
	%	\item[Phase 1] Constructing a simple polygon such that the convex rope problem is referred to as the problem of finding the shortest path joining $a$ and $b$ in the polygon;
	%	\item[Phase 2] Using three factors (f1) - (f3) of MMS  for the shortest path problem.
	%%	\item[Step 2] Partitioning of $\mathcal{D}$ into sub-polygons by \textit{cutting segments}.  At each cutting segment, we take  \textit{shooting points} to form the initial path;
	%%	\item[Step 3] 	We check a so-called \textit{collinear condition}  at shooting points to decide stop or not;
	%%	\item[Step 4]  If the algorithm does not stop, an update of shooting
	%%	points improves the paths. Then we obtain a set of new shooting points and a new path to go to Step 3.
	%\end{itemize}
	
	%The next part will provide a detailed exposition of  these phases in which Phase 1 is described in Sect.~\ref{sect:creating-polygon}  and Phase 2 is illustrated in Sect.~\ref{sect:partition},~\ref{sect:collinear_condition} and~\ref{sect:update}. 
	%
	%
	
	\subsection{Constructing a simple polygon $\mathcal{D}$}
	\label{sect:creating-polygon}
	% for convex rope problem ...

	%
	Take a rectangle $R$ that contains  $\mathcal{P}$ and has no edge touching $\mathcal{P}$. Since $b$ is visible from infinity, the ray $bb'$ always intersects with $R$ at only one point, say $c$.  The line segment $[b,c]$ partitions the non-simple region $R \setminus {\rm int}\mathcal{P}$ into a simple region, similar to that of~\cite{Toussaint1989} as follows. 
	%\begin{equation}
	%	\label{eq:construct-polygon}	\text{Insert $[b,c]$ and its copy which is }[\tilde{b},\tilde{c}] \text{ into the descriptions of the boundary of }R \setminus {\rm int}\mathcal{P}. \text{ A new simple polygon, denoted by } \mathcal{D}, \text{ is determined by a closed polyline including two parts as follows: one part, denoted by }B_1, \text{ is a polyline starting at } \tilde{b}, \text{ crossing } \tilde{c}, \text{ going along the boundary of } R \text{ with the direction of the arrow to come } c \text{ and going to } b; \text{ other, denoted by } B_2, \text{ is the polyline starting at } b \text{ moving along the boundary of } \mathcal{P} \text{ via clockwise order and returning to } \tilde{b} \text{ (see Fig.~\ref{fig:construct-polygon}(a)).}
	%	\end{equation}
	
	\begin{itemize}
		\item[(A)] We insert $[b,c]$ and its copy $[\tilde{b},\tilde{c}]$ into the descriptions of the boundary of $R \setminus {\rm int}\mathcal{P}$. Let $\mathcal{D}$ be a new simple polygon that is determined by a closed polyline including two parts: one part, denoted by $B_1$, is a polyline starting at $\tilde{b}$, crossing $\tilde{c}$, going along the boundary of $R$ with the direction of  arrows to come $c$ and going to $b$; other, denoted by $B_2$, is the polyline starting at $b$ going along the boundary of $\mathcal{P}$ via clockwise order and returning to  $\tilde{b}$. Thus $\mathcal{D}$ has the boundary $B_1 \cup B_2$ as shown in Fig.~\ref{fig:construct-polygon}(ii).
	\end{itemize}

	%	\begin{figure}[htp]
	%		\centering
	%		\begin{subfigure}[b]{0.5\textwidth}
	%			\centering
	%			\includegraphics[scale=0.28]{convert-to-MMS-new}
	%			\caption*{(i)}
	%		\end{subfigure}%
	%		\begin{subfigure}[b]{0.5\textwidth}
	%			\centering
	%			\includegraphics[scale=0.28]{is_visible_from}
	%			\caption*{(ii)}
	%		\end{subfigure}
	%		\caption{The counterclockwise convex rope starting $a$ and ending at $b$ is the shortest path joining $a$ and $b$ in the simple polygon $\mathcal{D}$, the case (i): $a$ is an extreme point of the convex hull of $\mathcal{P}$ and the case (ii): $a$ is visible from infinity.}
	%		\label{fig:construct-polygon}
	%	\end{figure}  

	The non-simple  region $R \setminus {\rm int}\mathcal{P}$ is converted into the simple polygon $\mathcal{D}$.
	By Remark~\ref{remark:a_is_on} given in Appendix, the convex rope starting at $a$ and ending at $b$ is indeed the shortest path joining $a$ and $b$ in $\mathcal{D}$. Therefore the convex rope problem is deduced to the shortest path problem in $\mathcal{D}$.
	We can also model the convex rope problem as the shortest path problem in another  polygon whose boundary includes  the convex hull of $\mathcal{P}$ instead of the rectangle $\mathcal{R}$, see~\cite{guibas-funnel}. However computing the convex hull of $\mathcal{P}$ needs linear time whereas the construction of rectangle $\mathcal{R}$ takes constant time when input points are within a given range.
	%	The next aim of this paper is to use MMS for solving the problem.
	According to Remark~\ref{remark:a_is_on} in Appendix, %in  Appendix (Sect.~\ref{app:detail-construct-polygon}), 
	from now on we focus on  finding the shortest path joining $\tilde{b}$ and $b$ in  $\mathcal{D}$.
	% Therefore the factors (f1)-(f3) are required for the multiple shooting. 

	\subsection{ The factor (f1): partitioning  the polygon $\mathcal{D}$} 
	\label{sect:partition}
	We split the polygon $\mathcal{D}$ into sub-polygons $\mathcal{D}_i$ by $N$
	cutting segments $\xi_1, \cdots,\xi_N$ ($ N \geq 1$) %(i.e., a shooting grid) 
	as follows
	\begin{equation}
		\label{eq:decomposition} %\tag{$P_{CVX}$}
		%\left\{
		\begin{array}{lll}
			\mbox{+ }\xi_i = [u_i, v_i] \subset \mathcal{D}  \mbox{ such that } ]u_i, v_i[ \subset {\rm int} \mathcal{D}, 
			\quad  u_i \in B_1 
			\mbox{ and } v_i \in B_2, \, \text{ for } 1\le i \le N; \\
			\mbox{+ } \xi_i \mbox{ strictly } \mbox{ separates } \tilde{b}
			\mbox{ and } b, \text{ for } 1\le i \le N;  \\
			\mbox{+ }  \xi_0 =\{\tilde{b}\}, \xi _{N+1} = \{b\}, [u_i, v_i] \, \cap \, [u_j, v_j] = \emptyset, \text{ for } i \neq j, \, 1\le i,j \le N;\\
			\mbox{+ } \mathcal{D}_i \mbox{ is bounded by } B_1,
			B_2, \mbox{ cut  segments }\xi_i \mbox{ and }\xi_{i+1}, \text{ for } 0\le i \le N;\\
			\mbox{+ } \mbox{int}\mathcal{D}_i \cap
			\mbox{int}\mathcal{D}_j = \emptyset, \mbox{  for } i \neq j,\, 0\le i,j \le N \mbox{  and } \mathcal{D}=\bigcup_{i=0}^N \mathcal{D}_i. 
		\end{array} %\right.
	\end{equation}
	
	% For simplicity, we may assume that the cutting segments $\xi_i$ are vertical line segments, for $i=1,2,\ldots, N$.
	%More details about an implementation of~(\ref{eq:decomposition})  will be seen in Sect.~\ref{sect:implementation}.
	%
	%By the condition $]u_i, v_i[ \subset {\rm int} \mathcal{D}, u_i \in B_2$ and $v_i \in B_1$, 
	
	%	A partition given by  (\ref{eq:decomposition}) is illustrated in Fig.~\ref{fig:cut-segments}(i). Clearly, $\mathcal{D}_i$ is  a simple polygon.
	A partition given by  (\ref{eq:decomposition}) is illustrated in Fig.~\ref{fig:MMS}(i). 
	%Clearly, $\mathcal{D}_i$ is  a simple polygon.	
	% In Fig.~\ref{fig:cut-segments}(ii), the decomposition is not a partition as in (\ref{eq:decomposition}), since $]u_i,v_i[$ and $]u_{i+1},v_{i+1}[$ do not lie completely in ${\rm int} \mathcal{D}$.
	%
	In the first step, 
	%we   initialize a set of ordered points by taking a point in each cutting segment.
	we take  a point in each cutting segment.
	Two consecutive initial points are connected by the shortest path 
	%joining these points in the polygon that indeed is the  shortest path joining them in corresponding sub-polygon. 
	in the corresponding sub-polygon. 
	%	The path received by combining these shortest  paths  is called the initial path of the algorithm. 
	The initial path is received by combining these shortest  paths.
	%	For convenience, these initial  points are chosen at the endpoints $v_i$ of cutting segments.
	For convenience, these initial  points are chosen at endpoints $v_i$ of $\xi_i$.
	For each iteration step, which is discussed carefully in next sections (\ref{sect:collinear_condition} and~\ref{sect:update}), we obtain an ordered set of points $\{a_i \vert \, a_i \in \xi_i,  i=1,2,\ldots,N\}$ and a path $\gamma =\cup_{i=0}^{N} \text{SP}(a_i, a_{i+1})$, where $\text{SP}(a_i, a_{i+1})$ is the shortest path joining $a_i$ and $a_{i+1}$ in $\mathcal{D}_i$ and $a_0 = \tilde{b}, a_{N+1} = b$. 
	Such a point $a_i$ is called a \textit{shooting point} and $\gamma$ is called \textit{the path formed by the set of  shooting points}
	% $\{a_i \vert \, i=1,2,\ldots,N\}$ 
	(see~\cite{An2017}).
	
	It is clear that the shortest path joining $a_i, a_{i+1}$ in $\mathcal{D}_i$ and that in $\mathcal{D}$ are identical. However, in the implementation, finding the shortest paths in these sub-polygons  makes more efficient use of time and memory than that in an entire polygon.  
	Here $\text{SP}(a_i, a_{i+1})$ could be computed by any known  algorithm for finding the shortest paths in simple polygons.
	Throughout this paper, when  we say  $\text{SP}(x, y)$ without further explanation, it means the shortest path joining $x$ and $y$ in the connected union of sub-polygons containing $x$ and $y$ of $\mathcal{D}$.
	
	%\begin{figure}[htp]
	%	\label{fig:cut-segments}
	%	\centering
	%		\begin{subfigure}[b]{0.4\textwidth}
	%	%	\centering
	%		\includegraphics[scale=0.5]{cut-segment1}
	%	\end{subfigure}%
	%	\begin{subfigure}[b]{0.4\textwidth}
	%%	\centering
	%	\includegraphics[scale=0.52]{cut-segment2}
	%	\end{subfigure}
	%	\caption{A partition of $\mathcal{D}$ given by (\ref{eq:decomposition}) as in (a), while the decomposition  as in (b) is not a partition given by (\ref{eq:decomposition}), since $]u_i,v_i[$ and $]u_{i+1},v_{i+1}[$ do not lie completely in ${\rm int} \mathcal{D}$.}
	%\end{figure}  
	
	%%%%%%%%%%%%%%%%%%%%%%%%
	%	\begin{figure}[htp]
	%		\centering
	%		\begin{subfigure}[b]{0.5\textwidth}
	%			\centering
	%			\includegraphics[scale=0.465]{cut-segment1}
	%			\caption*{(i)}
	%		\end{subfigure}%
	%		\begin{subfigure}[b]{0.5\textwidth}
	%			\centering
	%			\includegraphics[scale=0.48]{cut-segment2}
	%			\caption*{(ii)}
	%		\end{subfigure}
	%		\caption{A partition of $\mathcal{D}$ given by (\ref{eq:decomposition}) as in (i), while the decomposition  as in (ii) is not a partition given by (\ref{eq:decomposition}), since $]u_i,v_i[$ and $]u_{i+1},v_{i+1}[$ do not lie completely in ${\rm int} \mathcal{D}$ ($]u_i,v_i[$ and $]u_{i+1},v_{i+1}[$ contain some point(s) outside or on the boundary of $\mathcal{D}$).}
	%		\label{fig:cut-segments}
	%	\end{figure}  
	
	\subsection{The factor (f2): establishing and checking the collinear condition}
	\label{sect:collinear_condition}

	Firstly, we present a so-called  \textit{collinear condition} inspired by Corollary 4.3 in~\cite{HaiAnHuyen2019}. 
	% If the collinear condition  holds at all shooting points, the shortest path  is obtained from Proposition~\ref{prop:collinear-condition} in Sect.~\ref{sect:correctness}.
	% 	
	Recall that $a_i \in [u_i, v_i]$, for all $i=1,2,\ldots,N$, where $a_0 = \tilde{b}, a_{N+1} = b$.
	% Such points $a_i$ is called \textit{shooting points}. Let $\gamma =\cup_{i=0}^{N} \text{SP}(a_i, a_{i+1})$.
	%	For $i=1,2,\ldots,N$, 
	% Let $[s_1, a_i]$  	and $[a_i,s_2]$ be line segments of $\text{SP}(a_{i-1},a_i)$ and $\text{SP}(a_i,a_{i+1})$ that have the common endpoint $a_i$. 
	% 	The \textit{upper angle} created by $\text{SP}( a_{i-1}, a_i)$ and $\text{SP}(a_i, a_{i+1})$ with respect to $[u_i, v_i]$, denoted by  $\angle \left(  \text{SP}( a_{i-1}, a_i), \text{SP}(a_i, a_{i+1}) \right) $, is defined to be the total angle $ \angle s_1 a_i u'_i + \angle u'_i a_i s_2$,  where $u'_i$ is on the straight line $u_iv_i$ and $u_i$ is between $u'_i$ and $v_i$. 
	The \textit{upper angle} created by $\text{SP}( a_{i-1}, a_i)$ and $\text{SP}(a_i, a_{i+1})$ with respect to $[u_i, v_i]$, denoted by  $\angle \left(  \text{SP}( a_{i-1}, a_i), \text{SP}(a_i, a_{i+1}) \right) $, is defined to be the  angle of the polyline $\text{SP}( a_{i-1}, a_i) \cup \text{SP}(a_i, a_{i+1})$ at $a_i$ containing the ray $a_iu_i$.
	For each $i=1,2,\ldots,N$, the collinear condition states that:
	\begin{itemize}
		\item[(B)] If $a_i \in ]u_i, v_i[$, we say that $a_i$   satisfies the collinear condition when the upper angle created by $\text{SP}( a_{i-1}, a_i)$ and $\text{SP}(a_i, a_{i+1})$ is equal to $\pi$ (i.e, $\angle \left( \text{SP}( a_{i-1}, a_i), \text{SP}(a_i, a_{i+1})\right) =\pi$).
		\\
		If $a_i = v_i$, we say that $a_i$   satisfies  the collinear condition when the upper angle created by $\text{SP}( a_{i-1}, a_i)$ and $\text{SP}(a_i, a_{i+1})$ is at least $\pi$ (i.e, $\angle \left( \text{SP}( a_{i-1}, a_i), \text{SP}(a_i, a_{i+1})\right) \ge \pi$).
	\end{itemize}

	Let $\gamma^*$ be the shortest path joining $\tilde{b}$ and $b$ in $\mathcal{D}$, i.e., $\gamma^* = {\rm SP}(\tilde{b},b).$
	%	Note that, according to Remark~\ref{remark:a_is_on} in  Appendix (Sect.~\ref{app:detail-construct-polygon}), $\gamma^*$ never touches the boundary $B_1$ of $\mathcal{D}$ except  $b$. Thus Collinear Condition (B) does not include the case $a_i =u_i$.
	According to Remark~\ref{remark:a_is_on} in  Appendix,  Collinear Condition (B) does not include the case $a_i =u_i$.
	%	In the example given by  Fig.~\ref{fig:MMS}(i), Collinear Condition (B) holds at shooting points $a_{i-1}$ for the path  $\gamma^{current}$.

	\subsection{The factor (f3): updating  shooting points}
	\label{sect:update}

	%\begin{proposition}
	%	The sequence $<p_0,p_1,\ldots,p_{N+1}>$
	%	obtained by~(\ref{eq:formal-cvx-formulation}) is the unique solution of the
	%	shortest path problem.
	%\end{proposition}
	%
	Suppose that in $j^{th}$-iteration step, we have a set of shooting points $\{a_i\}_{i=1}^N$ and $\gamma =\cup_{i=0}^{N} \text{SP}(a_i, a_{i+1})$ is the path formed by the set of  shooting points, where $a_0 = \tilde{b}, a_{N+1} = b$. 
	If Collinear Condition (B) holds at all shooting points, by Proposition~\ref{prop:collinear-condition} in Sect.~\ref{sect:correctness}, the path obtained by the algorithm is the shortest, and our algorithm stops. Otherwise, we update shooting points to get a  path formed by the set of new shooting points in which  the sequence of paths obtained converges to  the shortest path, according to  Theorems~\ref{theo:global solution} and~\ref{theo:update_path}.
	Assuming that Collinear Condition (B) does not hold at all shooting points, we update  shooting points as follows
	
\begin{figure}[H]
		\centering
		\begin{subfigure}[r]{0.5\textwidth}
			\centering
			\includegraphics[scale=0.38]{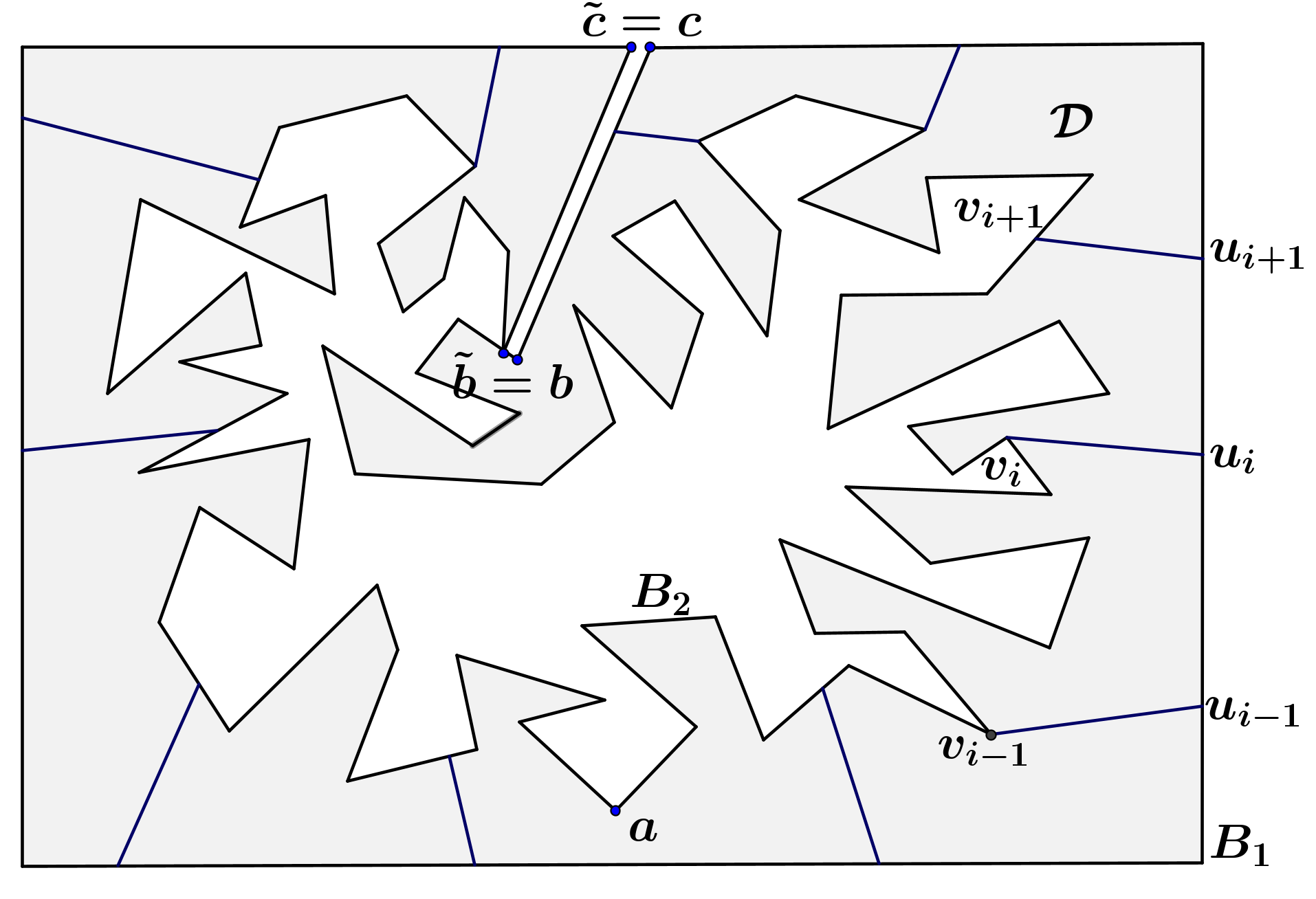}
			\caption*{(i)}
		\end{subfigure}%
		\begin{subfigure}[l]{0.5\textwidth}
			\centering
			\includegraphics[scale=0.38]{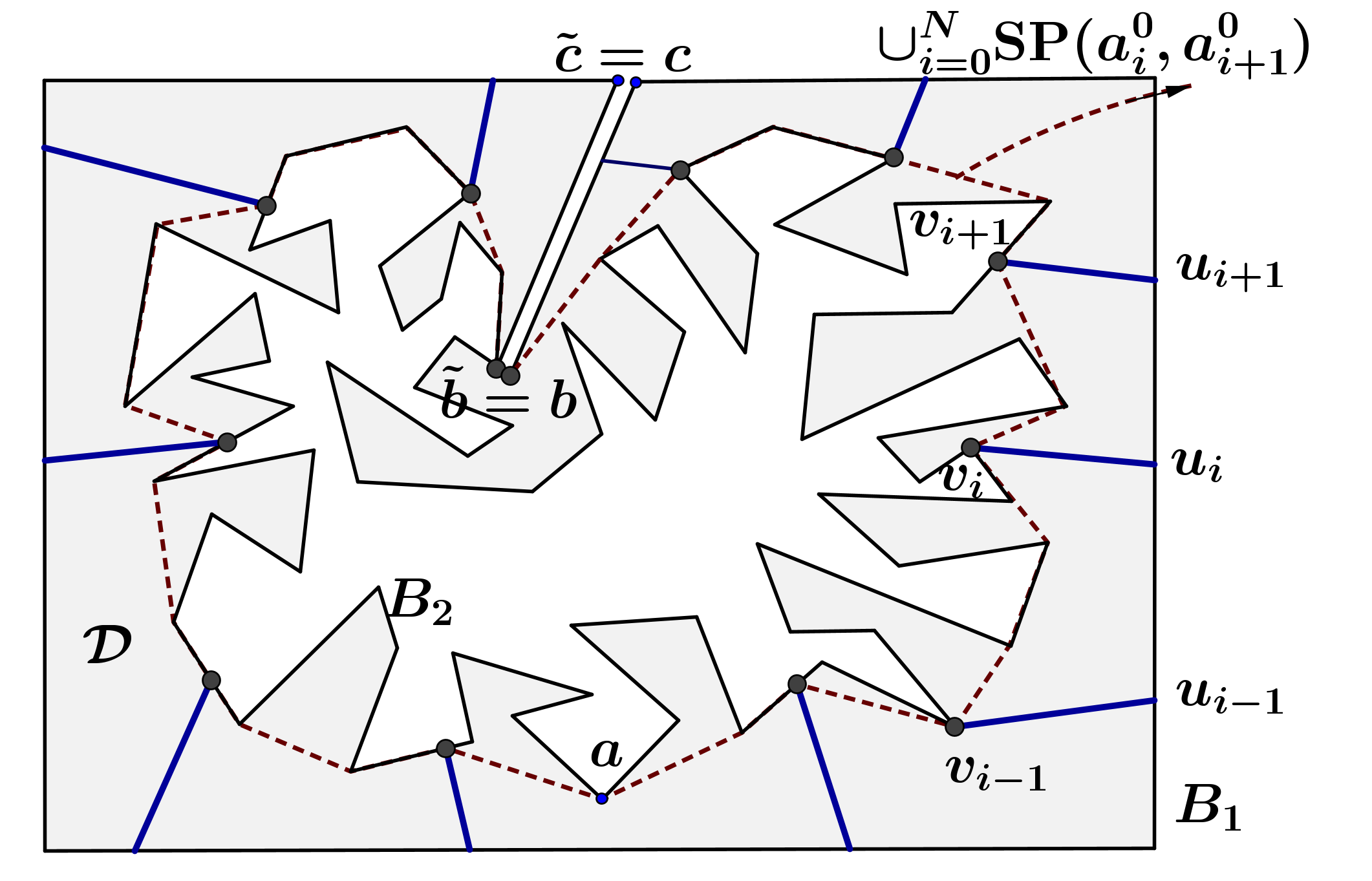}
			\caption*{(ii)}
		\end{subfigure}
		% 		\caption{(i): Illustration of partitioning $\mathcal{D}$ and (ii): The set of initial shooting points $\{a_i = v_i, i=0,1, \ldots, N+1\}$ (black big dots) and the path formed by  the set (dashed polyline).}
		% 	%	\label{fig:partition}
		% 	\end{figure}  
		% 	
		% 	
		% 	\begin{figure}[htp]
		
		\begin{subfigure}[l]{0.5\textwidth}
			\centering	\includegraphics[scale=0.75]{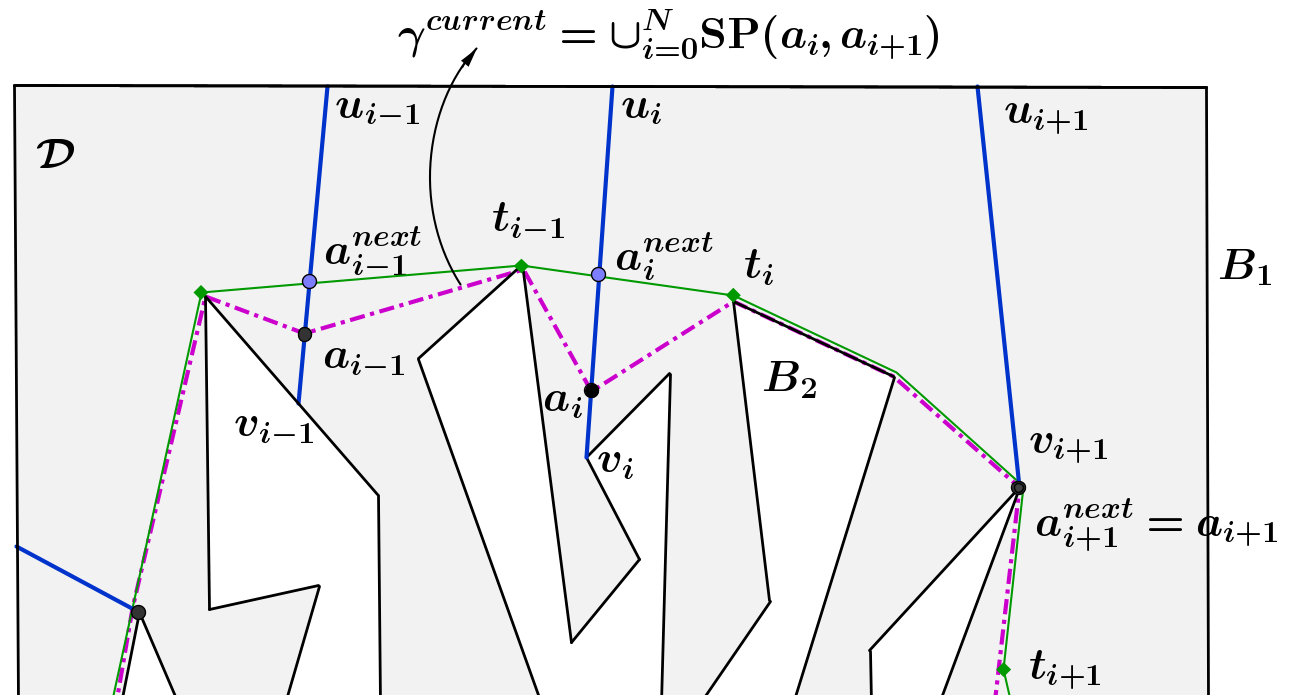}
			\caption*{(iii)}
		\end{subfigure}%
		\begin{subfigure}[r]{.5\textwidth}
			\centering
			\includegraphics[scale=0.75]{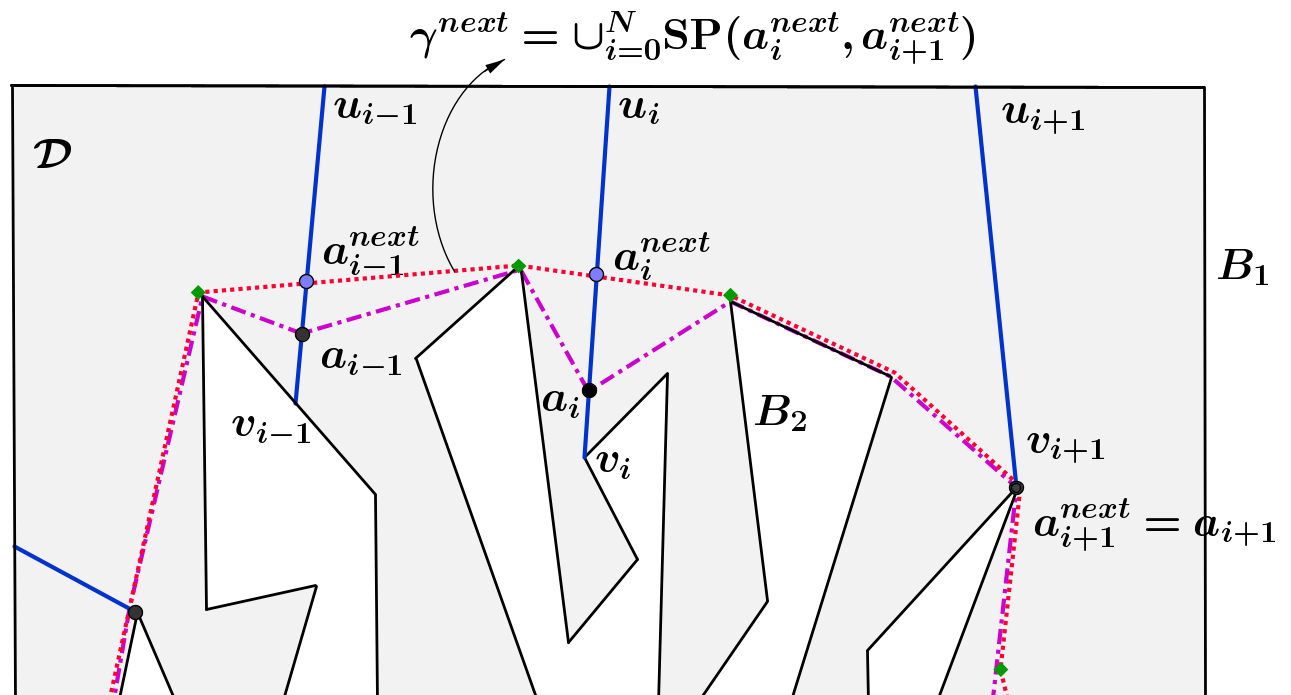}
			\caption*{(iv)}
		\end{subfigure}
		\caption{ (i): Illustration of partitioning $\mathcal{D}$; (ii): the set of initial shooting points consisting of $\tilde{b}$, $v_i, (i=1,2,\ldots,N)$ and $b$. (black big dots) and the path formed by  the set (dashed polyline); (iii): illustration of taking temporary points  and finding $a_i^{next}=\text{SP}(t_{i-1}, t_i) \cap \xi_i$; (iv): the update $\gamma^{current}$  (dashed polyline) to $\gamma^{next}$ (solid polyline).}
		\label{fig:MMS}
	\end{figure}

	%\begin{itemize}
	%	\item[(C)] 
	%	For $i=1,2 \ldots,N$, take $x_{i} \in \text{SP}(a_{i-1}, a_i), x'_i \in \text{SP}(a_i, a_{i+1})$ such that $x_i, x'_i$ are not identical to shooting points and $\gamma(x_{i-1}, x'_{i-1}) \setminus \{x_{i-1}, x'_{i-1}\}, \gamma(x_i, x'_i) \setminus \{x_i, x'_i\}$ and $\gamma(x_{i+1}, x'_{i+1}) \setminus \{x_{i+1}, x'_{i+1}\}$ are disjoint. Let $a'_i$ be a point of intersection of the shortest path $\text{SP}(x_i, x'_i)$ and $[u_i, v_i]$. The point of intersection exists uniquely since $[u_i, v_i]$ divides $\mathcal{D}$ into two parts, one that contains $x_i$ and one that contains $x'_i$. We update  shooting point $a_i$ to $a'_i$, see Fig.~\ref{fig:updating-way}(i).
	%\end{itemize}
	%For simplicity, we can choose $x'_{i-1} = x_i$ and $x'_i = x_{i+1}$. Therefore we rewrite the update of shooting points as follows:
	\begin{itemize}
		\item[(C)] 
		For $i=1,2 \ldots,N$, fix $t_{i-1} \in \text{SP}(a_{i-1}, a_i), t_i \in \text{SP}(a_i, a_{i+1})$ such that $t_{i-1}, t_i$ are not identical to shooting points. We call such point $t_{i-1}$ ($t_i$, resp.) the \textit{temporary point} with respect to $a_{i-1}$ ($a_i$, resp.). Let $a^{next}_i = \text{SP}(t_{i-1}, t_i) \cap [u_i, v_i]$. 
		%	The point of intersection exists uniquely since $[u_i, v_i]$ divides $\mathcal{D}$ into two parts, one that contains $t_{i-1}$ and one that contains $t_i$. We update  shooting point $a_i$ to $a^{next}_i$ (see Fig.~\ref{fig:updating-way}).
		
		The intersection point $a^{next}_i$ exists uniquely since $[u_i, v_i]$ divides $\mathcal{D}$ into two parts, one that contains $t_{i-1}$ and one that contains $t_i$. We update  shooting point $a_i$ to $a^{next}_i$ (see Fig.~\ref{fig:MMS}(iii)).
	\end{itemize}

	%	\begin{figure}[htp]
	%		\centering
	%		%	\begin{subfigure}[b]{0.5\textwidth}
	%		%		\centering
	%		%		\includegraphics[scale=0.4]{other_way_updating}
	%		%		\caption*{(i)}
	%		%	\end{subfigure}%
	%		%	\begin{subfigure}[b]{0.5\textwidth}
	%		%		\centering
	%		\includegraphics[scale=0.4]{other_way_updating}
	%		%	\caption*{(ii)}
	%		%	\end{subfigure}
	%		\caption{Illustration of taking $\{t_i, i=0,1, \ldots, N\}$   and finding $a_i^{next}$ that is the point of intersection of $\text{SP}(t_{i-1}, t_i)$ and $\xi_i$.}
	%		\label{fig:updating-way}
	%	\end{figure}  

	In particular, if $\text{SP}(a_i, a_{i+1})$ is a line segment, we choose $t_i$ as the midpoint of $\text{SP}(a_i, a_{i+1})$. Otherwise,  $t_i$ is chosen as an endpoint of $\text{SP}(a_i, a_{i+1})$ that is not identical to $a_i$ and $a_{i+1}$.
	We can update shooting points in which the collinear condition does not hold and keep the remaining shooting points. But the proposed algorithm will update all shooting points. Because if $a_{i}$ satisfies Collinear Condition (B), then the update gives $a_{i}^{next}=a_{i}$ due to 	 Proposition~\ref{prop:two-stop-condition} (in Sect.~\ref{sect:correctness}). 
	% For example in Fig.~\ref{fig:MMS}, Collinear Condition (B) holds at shooting points $a_{i-1}$ for $\gamma^{current}$.  	The update given by (C) also shows that $a_{i-1}^{next}=a_{i-1}$.
	\subsection{The Proposed Algorithm}

	%We next introduce the algorithm based multiple shooting approach with three factors established in Sect.~\ref{sect.MMS} to find $\gamma^*$ in $\mathcal{D}$. Note that we initialize $a_i^0 = v_i$, for $i=1,2,\ldots,N$, $a_0^0=\tilde{b}, a_{N+1}^0=b$.
	%\textit{\textbf{Main Algorithm:} Finding the convex rope joining $\tilde{b}$ and $b$ of the simple polygon $\mathcal{P}$}
	
	%	Our  algorithm is given as follows.
	%\begin{algorithm}[h!]
	\begin{algorithmic}[1]
		\Require A simple polygon $\mathcal{P}$, $b \in \mathcal{P}$  and its copy $\tilde{b}$, an integer number $N \ge 1$.
		\Ensure An approximate convex rope starting at $\tilde{b}$ and ending at $b$.
		% which does not enter the interior of $\mathcal{P}$.
		\State 
		Construct a simple polygon $\mathcal{D}$ with  the boundary $B_1 \cup B_2$  as shown in (A).
		\Comment{{\it \color{blue} see Fig.~\ref{fig:construct-polygon}}(ii)}
		\label{step1}
		\State Divide $\mathcal{D}$ into $\mathcal{D}_i$ satisfying~(\ref{eq:decomposition}) 
		by cutting segments $[u_i, v_i]$, where $u_i \in B_1, v_i \in B_2$ ($i=1,\ldots,N$)
		\Comment{{\it \color{blue} see Fig.~\ref{fig:MMS}(i)}}
		%\Comment{{\it \color{blue} partition}}
		\label{step2}
		\State Take an ordered set $\{a_i\}_{i=0}^{N+1}$ initial shooting points  consisting of $\tilde{b}$, $v_i, (i=1,2,\ldots,N)$ and $b$. Let $\gamma^{current}$ be the path formed by $\{a_i\}_{i=0}^{N+1}$.
		\Comment{{\it \color{blue} see Fig.~\ref{fig:MMS}(ii)}}
		\label{step3}
		\State $flag \leftarrow true$, $a_i^{next} \leftarrow a_i$ and $\gamma^{next}$ is the path formed by $\{a_i^{next}\}_{i=0}^{N+1}$.
		\label{step4}
		\State Call Procedure {\sc Collinear\_Update}($\gamma^{current}, \gamma^{next}, flag$)
		\Comment{{\it \color{blue} see Figs.~\ref{fig:MMS}(iii) and (iv)}}
		%\Comment{{\it \color{blue} given in Appendix (Sect.~\ref{sect:procedure})}}
		\label{step5}
		\State If Collinear Condition(B) holds at all shooting points, i.e., $flag= true$,		
		{\bf stop} and \Return $\gamma^{next}$. 
		\Comment{{\it \color{blue} the shortest path is obtained by Proposition~\ref{prop:collinear-condition}}}
		\label{step6}
		\State	 Otherwise, i.e, $flag = false$, then 		 
		$\gamma^{current} \leftarrow \gamma^{next}$ and  
		\textbf{go to } step~\ref{step5}.
		\label{step7}
		\State	\Return $\gamma^{next}$.
	\end{algorithmic}
	%\end{algorithm}
	For each iteration step, we need to check Collinear Condition (B) and update shooting points given by (C) to get a better path.  Then Procedure {\sc Collinear\_Update}($\gamma^{current}, \gamma^{next}, flag$) of the proposed algorithm performs checking condition (B) and updating the shooting point $a_i$ to $a_i^{next}$.
	%	 Figs.~\ref{fig:MMS}(iii) and (iv) illustrates updating shooting points given by (C). 
	%Theorem~\ref{theo:update_path} states the decreasing of  the sequence of lengths of paths obtained by the procedure.

	%
	\begin{algorithm}[htp]
		\begin{algorithmic}[1]
			\Procedure{Collinear\_Update}{$\gamma^{current}, \gamma^{next}, flag$} \label{alg:update}
			%
			%		\Require $\gamma^{current}=\cup_{i=0}^N \text{SP}(a_i,a_{i+1})$, $flag=false$.
			%		\Comment{{\it \color{blue} $flag=false$ then Collinear Condition (B) does not hold at some index $i, (1 \le i \le N)$}}
			%		\Ensure Determined if   $flag=true$  or $flag=false$  and the path formed by the set of new shooting points $\gamma^{next}=\cup_{i=0}^N \text{SP}(a_i^{next},a_{i+1}^{next})$.
			%\Statex
			\State  $flag \leftarrow true$
			\For {$i=0,1, \ldots, N$}
			\State
			\label{step:temp_points} Take temporary points $t_i$ on  $\text{SP}(a_i,a_{i+1})$ such that \par
			\quad	if $\text{SP}(a_i, a_{i+1})$ is a line segment, $t_i$ is the midpoint of $\text{SP}(a_i, a_{i+1})$, \par
			\quad	for if not,  $t_i$ is  an endpoint  of $\text{SP}(a_i, a_{i+1})$ that is not identical to $a_i$ and $a_{i+1}$.
			\If { $a_i \in ]u_i,v_i[$} 
			\If {$\angle \left( \text{SP}( a_{i-1}, a_i), \text{SP}(a_i, a_{i+1}) \right) =\pi$} 
			\Comment{{\it \color{blue}  Collinear Condition {\rm (B)} is checked}}
			\State Set $a_i^{next} \leftarrow a_i$
			\Else \ call $\text{SP}(t_{i-1}, t_i)$ 
			\label{step-prc8}
			\State Set $a_i^{next}$ be the intersection point of  $\text{SP}(t_{i-1}, t_i)$ and $[u_i, v_i]$
			\Comment{{\it \color{blue} due to {\rm (C)}}}
			\State $flag \leftarrow false$
			\EndIf 
			\Else \  
			\Comment{{\it \color{blue} $a_i =v_i$}}
			\If {$\angle \left( \text{SP}( a_{i-1}, a_i), \text{SP}(a_i, a_{i+1}) \right) \ge \pi$} 
			\Comment{{\it \color{blue}  Collinear Condition {\rm (B)} is checked}}
			\State Set $a_i^{next} \leftarrow a_i$
			\Else \ call $\text{SP}(t_{i-1}, t_i)$ 
			\label{step-prc14}
			\State Set $a_i^{next}$ the intersection point of  $\text{SP}(t_{i-1}, t_i)$ and $[u_i, v_i]$
			\Comment{{\it \color{blue} due to {\rm (C)}}}
			\State $flag \leftarrow false$
			\EndIf 
			\EndIf
			\EndFor
			\State $a_0^{next}=\tilde{b}, a_{N+1}^{next}=b$
			\EndProcedure
		\end{algorithmic}
	\end{algorithm}

	\section{The Correctness of the Proposed Algorithm}
	\label{sect:correctness}
	%Let $\{a_i, i=1,2,\ldots,N\}$ be a set of shooting points, where $a_0 = \tilde{b}, a_{N+1} = b$. Let $\gamma =\cup_{i=0}^{N} \text{SP}(a_i,a_{i+1})$.
	%	In this section, we introduce some results in which their proof are given in  Appendix.
%	Due to the properties of the shortest path joining two points in a simple polygon as shown in  Appendix, we have
%	\begin{remark}
%		\label{remark: optimal_solution}	
%		Suppose that the shortest path $\gamma^*$  intersects with $\xi_i$ at $a_i$, for $i=1,2, \ldots,N$.  Then Collinear Condition (B) holds at   $a_i$,  for all $i=1,2,\ldots,N$.
%	\end{remark}
	
	%	According to Remark~\ref{remark: optimal_solution}, the shortest path $\gamma^*$ satisfies Collinear Condition (B). 
	The following proposition shows  that the path satisfying Collinear Condition (B) is the shortest path $\gamma^*$ joining $\tilde{b}$ and $b$. 
	
	\begin{proposition}
		\label{prop:collinear-condition}
		Let $\gamma^{current}$ be a path formed by a set of shooting points, which joins $\tilde{b}$ and $b$ in $\mathcal{D}$.
		If Collinear Condition (B) holds at all  shooting points, then $\gamma^{current} = \gamma^*.$
	\end{proposition}
	
	%	The proof of Proposition~\ref{prop:collinear-condition} is given in  Appendix.

	\begin{proposition}
		\label{prop:two-stop-condition}
		Collinear Condition (B) holds  at  $a_i$ if and only if ${\rm SP}(t_{i-1}, t_i)\cap [u_i, v_i] = a_i$, where $t_{i-1}$ and $t_i$ are temporary points w.r.t. $a_{i-1}$ and $a_i$.  
	\end{proposition}

	\begin{theorem}
		\label{theo:global solution}
		The algorithm is convergent. It means that if the algorithm stops after a finite number of iteration steps (i.e.,  Collinear Condition (B) holds at all shooting points), we obtain the shortest path, otherwise, 
		the sequence of paths $\{\gamma^{j}\}$ converges to $\gamma^*$ in Hausdorff distance (i..e., $d_{\mathcal{H}}(\gamma^{j}, \gamma^*) \rightarrow 0$ as $j \rightarrow \infty$, where $\gamma^{j}$ be a path obtained by the  algorithm in $j^{th}$-iteration step). 
		%Therefore $\lim\limits_{j\rightarrow \infty} \gamma^j = \gamma^*$.
	\end{theorem}
	%	The proof of Theorem~\ref{theo:global solution} is given in  Appendix.
	%, where $\lim\limits_{j\rightarrow \infty} \gamma^j$ is the closed limit of a sequence of paths  in $\mathbb{R}^2$ (see Sect. 4.3 \cite{Papadopoulos2005}). 
	Since the  function of the length of paths is lower semi-continuous but not continuous, the event $d_{\mathcal{H}}(\gamma^{j}, \gamma^*) \rightarrow 0$ as $j \rightarrow \infty$ does not ensure the  convergence of the sequence of lengths of these paths. However, this holds for the sequence of paths  obtained by the proposed algorithm.  This is shown in the below theorem.
	
	% 	\begin{proposition}
	% 		\label{prop:update_path} The update of the propose algorithm by the factor (f3), i.e.,  Procedure \break {\sc Collinear\_Update}($\gamma^{current}, \gamma^{next}, flag$), gives $	l\left(\gamma^{next} \right) \le l\left( \gamma^{current}\right)$.
	%% 		\begin{align*}
	%% 			%	\label{eq:update_length_descreasing}
	%% 			l\left(\gamma^{next} \right) \le l\left( \gamma^{current}\right)
	%% 		\end{align*}
	% 		%
	% 		%Let $\{a_i, i=1,2,\ldots,N \}$ et of shooting points, where $a_0 = a, a_{N+1} = b$. Let $\gamma =\cup_{i=0}^{N} \text{SP}(a_i,a_{i+1})$.	For $i=1,2, \ldots,N$ if the collinear condition (B) is not satisfied at $a_i$, we update $a_i$ to new shooting point $a'_i$ by (B) in Sect.~\ref{sect:update}. Otherwise, the collinear condition (B) is satisfied at $a_i$, we denote $a'_i = a_i$. Let $\gamma' =\cup_{i=0}^{N} \text{SP}(a'_i,a'_{i+1})$, where $a'_0=c, a'_{N+1}=b$. Then $l(\gamma') \le l(\gamma).$
	% 		%
	% 		Thus the  sequence of lengths of paths  obtained by our algorithm is convergent. 	If Collinear Condition (B) does not hold at some shooting point, then the inequality
	% 		above is strict.
	% 	\end{proposition}
	
	\begin{theorem}
		\label{theo:update_path} The  sequence of lengths of paths  obtained by the proposed algorithm has the following properties:
		\begin{itemize}
			\item[(a)] 
			The  sequence is decreasing, i.e.,  Procedure  {\sc Collinear\_Update}($\gamma^{current}, \gamma^{next}, flag$) gives $	l\left(\gamma^{next} \right) \le l\left( \gamma^{current}\right)$. If Collinear Condition (B) does not hold at some shooting point, then the inequality
			above is strict.
			\item[(b)]	The  sequence converges to  $l(\gamma^*)$. 	
		\end{itemize}
	\end{theorem}
	
	%	The proof of Theorem~\ref{theo:update_path} is given in  Appendix.
	%  According to Proposition~\ref{prop:collinear-condition} and Theorem~\ref{theo:update_path}, we obtain
	% \begin{corollary}
	% 	\label{coro:convergent}
	% 	If Collinear Condition (B) holds, the algorithm gives the shortest path joining $\tilde{b}$ and $b$ in $\mathcal{D}$. Otherwise, the sequence of  lengths of  paths obtained by the algorithm is convergent.
	% \end{corollary}
	% \begin{proof}
	% 	When the algorithm performs, if Collinear Condition (B) holds, due to Proposition~\ref{prop:collinear-condition}, the shortest path joining $\tilde{b}$ and $b$ in $\mathcal{D}$ is obtained. Otherwise, by Theorem~\ref{theo:update_path}, the sequence of  lengths of  paths  found by the algorithm   decreases strictly. Then the sequence of lengths  of paths is convergent.
	% \end{proof}

	\begin{theorem}
		\label{theo:complexity}
		Given a fixed number $N$ of cutting segments, the proposed algorithm runs in $O(kn)$ time, where $n$ is the number of vertices of $\mathcal{P}$, $k$ is the  number of iterations to get the required path joining $\tilde{b}$ and $b$.
	\end{theorem}
	
	The proofs of the above propositions and theorems  are given in  Appendix.

	\section{Numerical Experiments}
	\label{sect:implementation}
	In the implementation, assume that the boundary of $\mathcal{P}$ consists of two polylines with integer coordinates that are monotone w.r.t. $x$-coordinate axis, and the cutting segments $\xi_i$ are parallel to $y$-coordinate axis.
	The monotone condition appears because we use a procedure for finding the shortest path joining two consecutive shooting points  inside sub-polygons $\mathcal{D}_i$ 
	as shown in~\cite{AnHoai2010}. 
	%Furthermore our code runs in the case coordinates of vertices of $\mathcal{P}$ are integer numbers.
	%
	%The part of  the boundary of convex hull of $\mathcal{P}$ from $a$ to $a$  is computed by using Melkman's algorithm~\cite{melkman}.
	%The initialization of shooting points, noted $a_i^0 = v_i, i=1,2,\ldots,N$. 
	%
	The algorithm is implemented in C++ programming language. The codes are compiled and executed 
	by GNU Compiler Collection under platform Ubuntu Linux 18.04, Processor 2.50GHz Core i5. 
	% We use a parameter $\epsilon =1e-6$
	%For the convenience of the implementation, we use cutting segments which are parallel to $y$-axis. According to 
	%
	
\begin{figure}
		\centering
		\begin{subfigure}{.5\textwidth}
			\centering
			\includegraphics[scale=0.18]{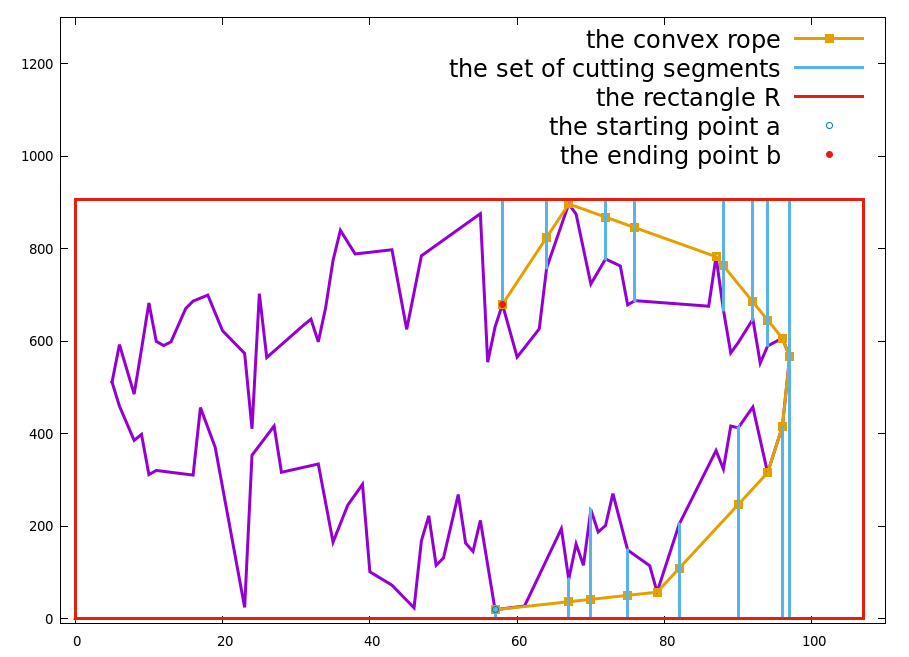}
			\caption*{(i)}
		\end{subfigure}%
		\begin{subfigure}{.5\textwidth}
			\centering
			\includegraphics[scale=0.18]{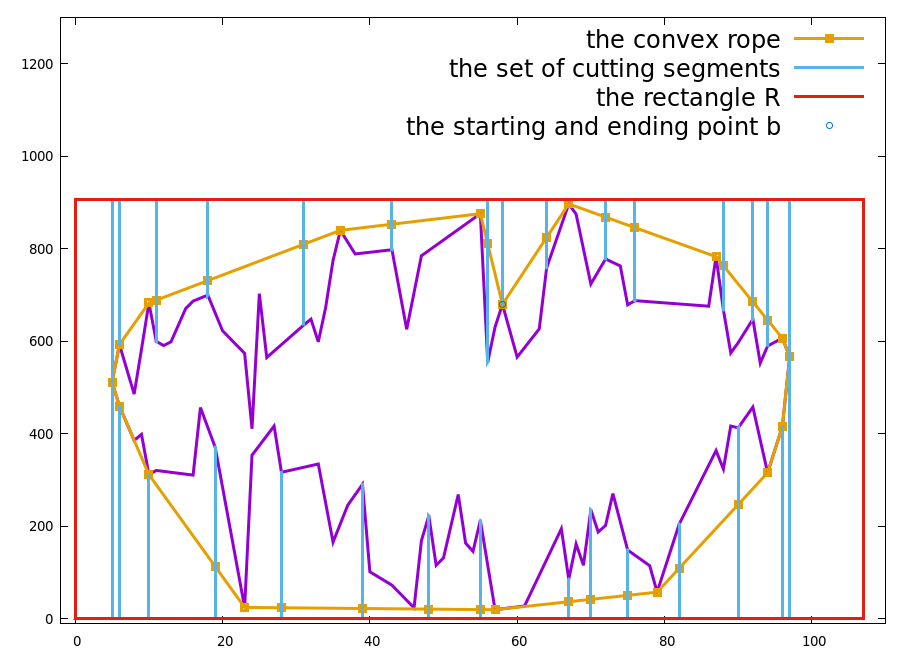}
			\caption*{(ii)}
		\end{subfigure}
		\caption{Convex ropes of a polygon with 200 vertices:
			(i) the convex rope starting at a vertex of the convex hull of $\mathcal{P}$ and ending at a point that is visible from infinity; (ii) the convex rope starting  and ending at the same point that is visible from infinity.}
		\label{fig:implemetation2}
	\end{figure}

\begin{figure}[htp]
		%	\label{fig:implemetation2}
		%	\centering
		\begin{subfigure}{.5\textwidth}
			\centering
			\includegraphics[scale=0.172]{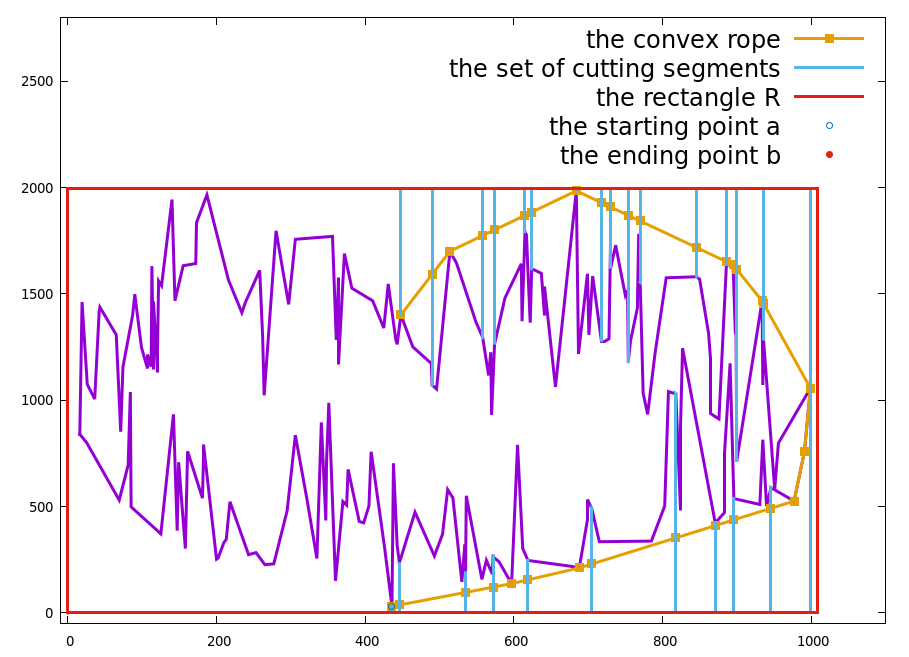}
			\caption*{(i)}
		\end{subfigure}%
		\begin{subfigure}{.5\textwidth}
			\centering
			\includegraphics[scale=0.18]{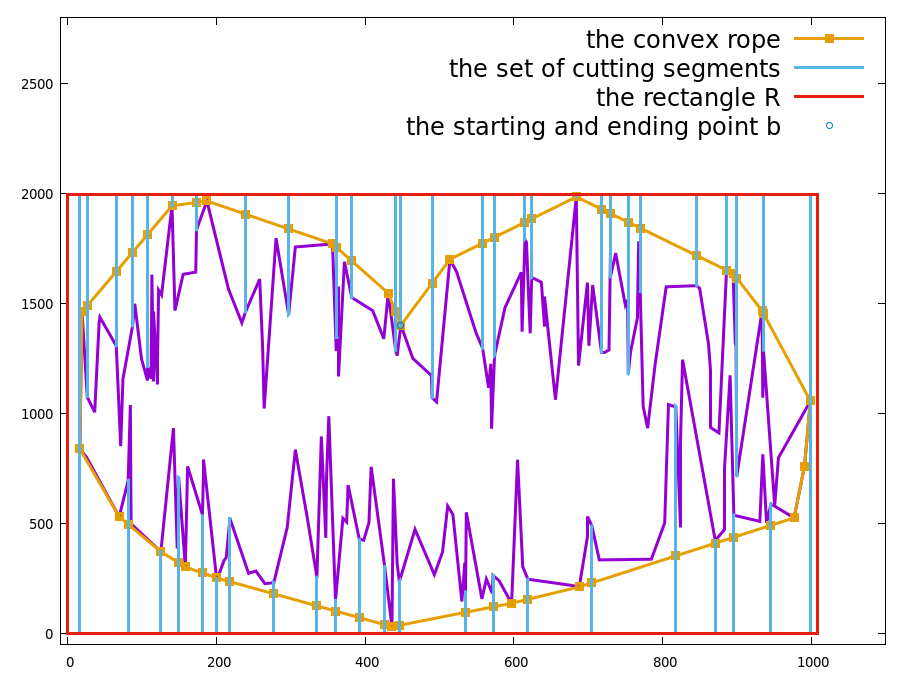}
			\caption*{(ii)}
		\end{subfigure}
		\caption{Convex ropes of a polygon with 400 vertices: (i) the convex rope starting at a vertex of the convex hull of $\mathcal{P}$ and ending at a point that is visible from infinity; (ii) the convex rope starting  and ending at the same point that is visible from infinity.}
		\label{fig:implemetation3}
	\end{figure}

\begin{figure}[H]
		%\label{fig:implemetation3}
		%	\centering
		\begin{subfigure}{.5\textwidth}
			\centering
			\includegraphics[scale=0.176]{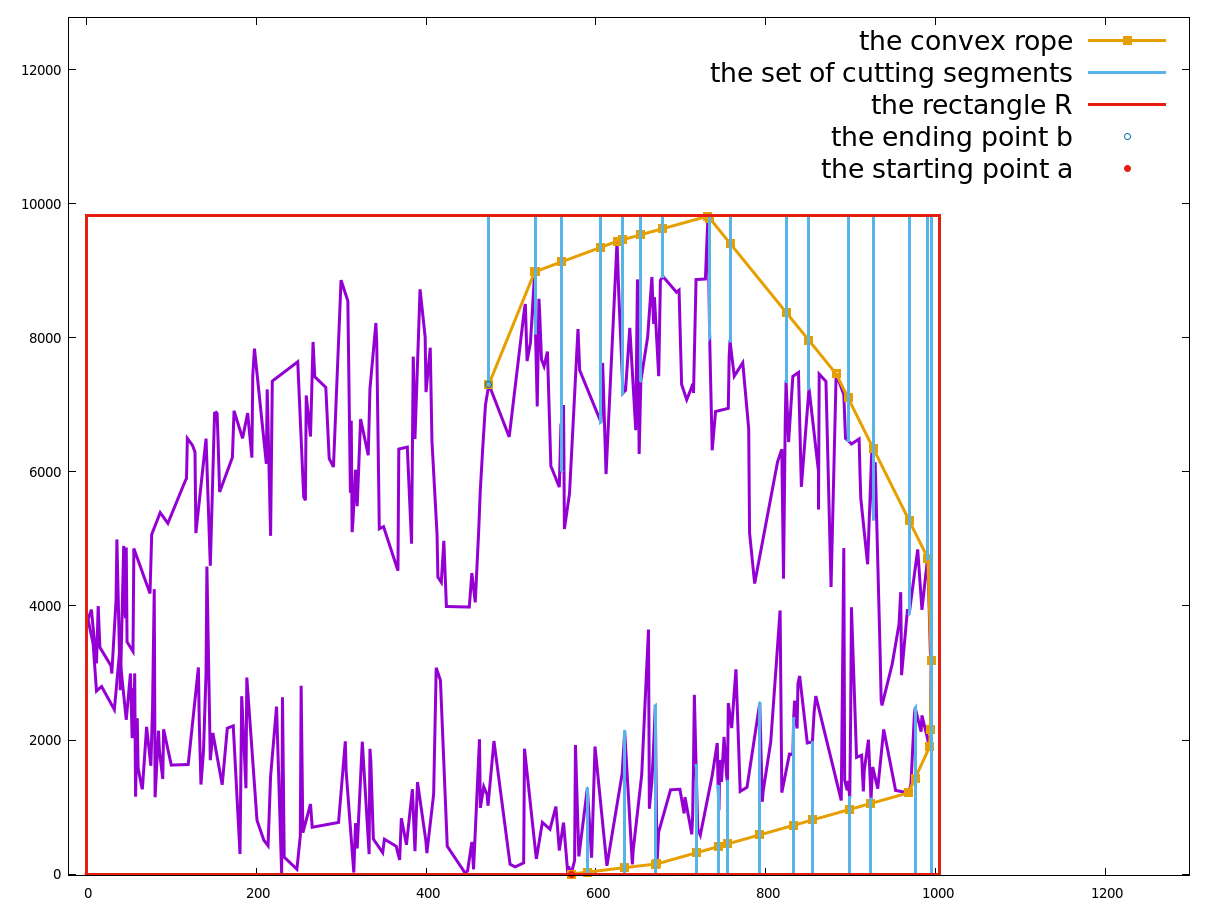}
			\caption*{(i)}
		\end{subfigure}%
		\begin{subfigure}{.5\textwidth}
			\centering
			\includegraphics[scale=0.176]{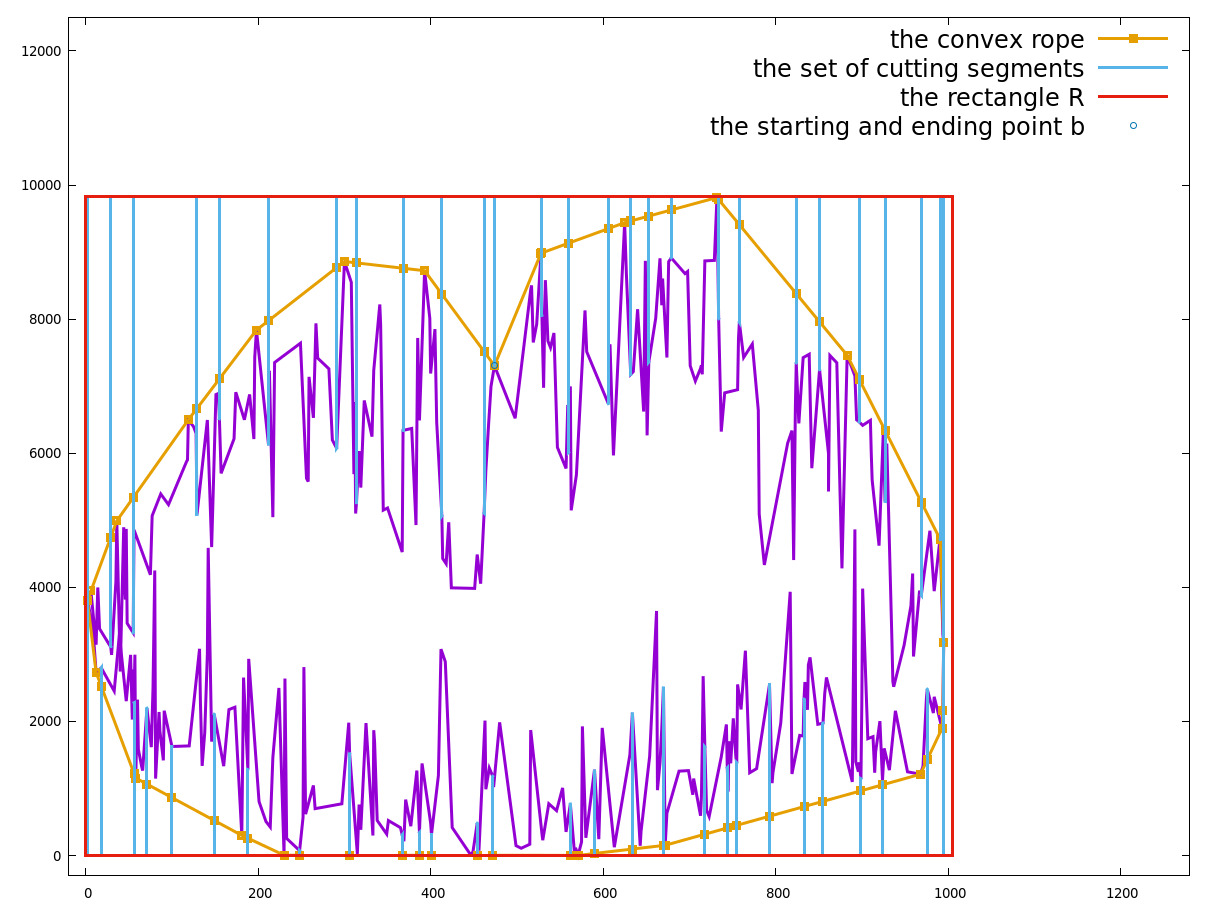}
			\caption*{(ii)}
		\end{subfigure}
		%\label{fig:implemetation4}
		\caption{Convex ropes of a polygon with 800 vertices: (i) the convex rope starting at a vertex of the convex hull of $\mathcal{P}$ and ending at a point that is visible from infinity; (ii) the convex rope starting  and ending at the same point that is visible from infinity.}
		\label{fig:implemetation4}
	\end{figure}

\begin{SCfigure}
		%
		%	\centering
		\includegraphics[width=0.56\textwidth]{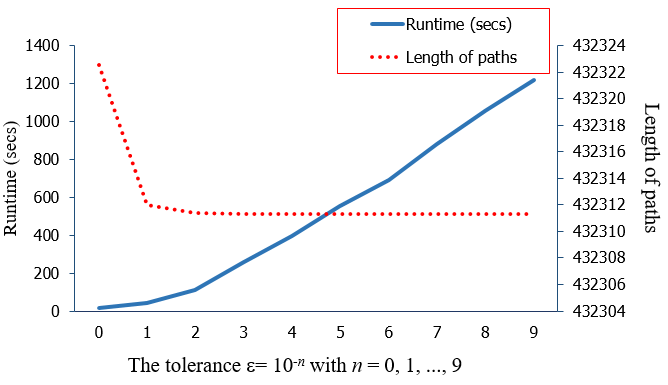}
		\caption{The dependence of the runtime of the proposed algorithm and the length of corresponding convex ropes on $\epsilon$, where the convex ropes  of a polygon having 3000 vertices starting  and ending at the same point that is  visible from infinity.}
		\label{fig:effect-tolerance}
		%	\end{figure}
	\end{SCfigure}

	To check whether Collinear Condition (B) holds or not at a shooting point $a_i$, we calculate the upper angle at $a_i$ or determine if $a_i^{next}$ coincides with $a_i$ or not. Two ways are the same by Proposition~\ref{prop:two-stop-condition}. Here we use the latter way.
	We need a tolerance, say $\epsilon$, to check for the coincidence of points when implementing. Collinear Condition (B) holds at $a_i$  if $ \Vert a_i^{next} - a_i \Vert < \epsilon$, for all $1 \le i \le N$.
	Moreover, by Lemma~\ref{lem:HaiAn2011}  in Appendix, the condition $ \Vert a_i^{next} - a_i \Vert < \epsilon$, $1 \le i \le N$ also applies to determine $\{\gamma^j\}$ converging to $\gamma^*$. 
	The actual runtime  of our algorithm grows with the decreasing of $\epsilon$. This effect is given in Table~\ref{tabl:tolerance} and Fig~\ref{fig:effect-tolerance}, where $\epsilon$ changes from $10^0$ to $10^{-9}$, and the problem is to find an approximate convex rope starting at the copy $\tilde{b}$ of a fixed point $b$ and ends at $b$ of the polygon having 3000 vertices and 200 cutting segments.
	Figs.~\ref{fig:implemetation2},~\ref{fig:implemetation3}, and~\ref{fig:implemetation4} show the results when $\epsilon = 10^6$ and  the number of vertices of $\mathcal{P}$ is  $100$, $400$, and $800$, respectively. Herein, the convex ropes start at $a$ and end at $b$, where $a$ is a visible infinity point that plays the same role as $b$ ($a= \tilde{b}$) or a vertex of the convex hull of $\mathcal{P}$.
	
	%In  the case is $ \epsilon = 10^{-6}$, Figures at page~\pageref{fig:implemetation1} illustrate convex ropes starting at $a$ and ending at $b$, where $a$ is visible infinity point, that plays the same role as $b$ ($a= \tilde{b}$) or $a$ is an extreme point of the convex hull of $\mathcal{P}$, where the number of vertices of $\mathcal{P}$ is  $100$, $400$ and $800$, respectively.

	%Table~\ref{tabl:tolerance} lists that when $\epsilon = 10^{-10}$, Collinear Condition (B) does not hold. Then we use a tolerance, say $\epsilon_h$, to ensure that the code terminates if $ | l(\gamma^{next})-l(\gamma^{current}) | < \epsilon_h$. The results in Table~\ref{tabl:tolerance} corresponds to $\epsilon_h = 10^{-4}$.
	
	\newpage

	
	\section{Conclusion}
	\label{sect:conclusion}
	We presented the use of MMS with three factors (f1)-(f3) for approximately solving the convex rope problem.
	%The principal advantages of the algorithm for finding the convex rope can be summarized as follows. 
	%	Although the ``oracle'' condition does not hold for polygon $\mathcal{D}$, unlike results in~\cite{AnHaiHoai2013}, we prove that the shortest path joining $\tilde{b}$ and $b$ inside $\mathcal{D}$ can be obtained by MMS. 
	Although the ``oracle'' condition does not hold for polygon $\mathcal{D}$ as in~\cite{AnHaiHoai2013},  we dealt with three open questions presented in Sect.~\ref{sec:intro}.
	% and obtain fully theoretical results.
	%Using MMS with slight modifications on determining  the collinear condition and updating new paths can obtain  a sequence of paths which converges to the  optimal solution by the Hausdorff distance. In addition, if the colinear condition holds at all shooting points, then the shortest path is well determined. Otherwise, the sequence of lengths of paths found by the algorithm is strictly decreasing and convergent.
	%We presented the use of the multiple shooting method for approximately solving the convex rope problem. In this work, slight modifications of the collinear condition and shooting point updating were adapted for the problem so that a sequence of paths obtained by the algorithm converges to the shortest path by means of the Hausdorff distance.
	%
	%
	We  proved that the shortest path is well determined if the collinear condition holds at all shooting points. Otherwise, the sequence of paths obtained from new shooting points by the update of the method converges in Hausdorff distance to the  shortest path. 
	The proposed algorithm was implemented in C++ and some numerical examples were given to show that the method is suitable for the problem. On the other side, as an advantage of MMS by the factor (f1)~\cite{AnHaiHoaiTrang2014}, instead of solving the problem for the whole formed polygon, we deal with a set of smaller partitioned subpolygons. 
	Consequently, the memory consuming of the algorithm should be low. This is useful when considering to deploy the method in robotic applications with limited computing resources. This will be a further research in the future. 
	
	%\begin{table}[htp]
	%	\centering
	%	\caption{ The dependence of the runtime of the algorithm on the tolerances}
	%	\label{tabl:tolerance}
	%	\begin{tabular}[b]{r|r|r}
	%		\toprule    %[0.3pt] 
	%		Tolerance $\epsilon$ &  The runtime (in sec.)  & The runtime (in sec.) \\
	%		 &  without $\epsilon_h$  & with  $\epsilon_h =10^{-4}$   \\
	%		\midrule    %[0.3pt]
	%	%	\noalign{\smallskip}
	%		$10^0$ & $14.26$ & $18.52$ \\
	%		$10^{-1}$ & $44.33$ & $39.50$ \\
	%		$10^{-2}$ & $106.96$ & $87.09$ \\
	%		$10^{-3}$ & $257.45$ & $87.10$\\
	%		$10^{-4}$ & $398.17$ & $87.03$ \\
	%		$10^{-5}$ & $553.05$ & $84.94$\\
	%		$10^{-6}$ & $693.02$ & $85.89$ \\
	%		$10^{-7}$ & $881.29$ & $87.93$\\
	%		$10^{-8}$ & $1053.78$ & $86.22$ \\
	%		$10^{-9}$ & $1217.04$ & $91.36$ \\
	%		$10^{-10}$ & Collinear Cond. does not hold & $88.70$ \\
	%		\bottomrule %[0.3pt]
	%	\end{tabular}
	%\end{table}
	\begin{SCtable}	
		%	\begin{table}[htp]
		\centering
		\caption{ The actual runtime of the proposed algorithm and the length of corresponding convex ropes when $\epsilon$ changes, where the convex ropes of a polygon having 3000 vertices start  and end at the same point that is  visible from infinity.}
		\label{tabl:tolerance}
		\begin{tabular}[b]{r|r|r}
			\toprule    %[0.3pt] 
			Tolerance $\epsilon$ &  The runtime (secs)   & Length of paths  \\
			\midrule    %[0.3pt]
			%	\noalign{\smallskip}
			$10^0$ & $14.26$ & $432322.554$ \\
			$10^{-1}$ & $44.33$ & $432312.019$ \\
			$10^{-2}$ & $106.96$ & $432311.361$ \\
			$10^{-3}$ & $257.45$ & $432311.340$ \\
			$10^{-4}$ & $398.17$ & $432311.339$ \\
			$10^{-5}$ & $553.05$ & $432311.339$ \\
			$10^{-6}$ & $693.02$ & $432311.339$ \\
			$10^{-7}$ & $881.29$ & $432311.339$ \\
			$10^{-8}$ & $1053.78$ & $432311.339$ \\
			$10^{-9}$ & $1217.04$ & $432311.339$ \\
			\bottomrule %[0.3pt]
		\end{tabular}
		%	\end{table}
	\end{SCtable}	
	
	%	\begin{figure}[H]
	
	%\newpage

	\section*{Acknowledgment}
	This research is funded by Vietnam National
	University Ho Chi Minh City (VNU-HCM) under grant number T2022 - 20 -01.
	
	The first and third authors acknowledge Ho Chi Minh City University of Technology (HCMUT), VNU-HCM for supporting this study.
	
The second author was funded by Vingroup Joint Stock Company and supported by the Master/ PhD Scholarship Programme of Vingroup Innovation Foundation (VINIF), Vingroup Big Data
Institute (VinBigdata), code VINIF.2021.TS.068.

	% (a) The triangulation of the polygon is not neces-
	%sary. (b) The location of vertices of the polygon from a toin clockwise order does not play any role in the algorithm.
	%Moreover, computing the convex hull of the polygon P is
	%not needed in some cases of a and b (such a and b are given
	%in Section 4).
	
	% references
	
	%
	\appendix
	\section*{Appendix}
	\renewcommand{\thesubsection}{\Alph{subsection}}

	\subsection{Shortest paths and convex ropes} 
	\label{app:detail-construct-polygon}
	With the construction of $\mathcal{D}$ as shown by (A) in Sect.~\ref{sect:creating-polygon},  vertices of $\mathcal{D}$ on $B_1$  except for $b$ and $\tilde{b}$ are convex (see Fig.~\ref{fig:construct-polygon}(ii)). 
	% 	As $\text{SP}(a,b)$ in the polygon $\mathcal{D}$ is a polyline whose vertices except for $a$ and $b$ are reflex, $\text{SP}(a,b)$ does not pass through any vertices of $B_1$. Hence we 
	We have the following
	\begin{remark}
		\label{remark:boundaryB_2}
		The shortest path joining $a$ and $b$ in $\mathcal{D}$ does not depend on the rectangle $R$, i.e, we can replace $R$ with any  rectangle that contains $\mathcal{P}$ and has no edge touching $\mathcal{P}$ to obtain another polygon $\mathcal{D}'$ in which shortest paths joining $a$ and $b$ in $\mathcal{D}$ and $\mathcal{D}'$ are the same. Furthermore, $\text{SP}(a,b)$ never touches  $B_1$ of $\mathcal{D}$ except for $b$.
	\end{remark}
	%  Le 6/1/2022
	%	Since a convex rope is the shortest path joining two points in $\mathcal{D}$ that is unique, a convex rope starting at $a$ and ending at $b$ exists uniquely.
	%	Recall that $\gamma^*$ is the shortest path joining $\tilde{b}$ and $b$ in $\mathcal{D}$.

	% 	\medskip
	%	Since $a$ is an extreme point of the convex hull of $\mathcal{P}$, we have the following
	
	\begin{remark} \label{remark:a_is_on} 
		The convex rope starting at $a$ and ending at $b$ of $\mathcal{P}$ can be obtained by solving directly the problem of finding the shortest path joining $a$ and $b$ inside  $\mathcal{D}$.	
		If $a$ is a vertex of the convex hull of $\mathcal{P}$, then the shortest path contains $a$ and never touches the boundary $B_1$ of $\mathcal{D}$ except for $\tilde{b}$ and $b$. More generally, if $a$ is visible from infinity, then $\tilde{b}$ plays the same role as $a$ (see Fig.~\ref{fig:construct-polygon}(ii)). Therefore  from now on, we only focus on finding the shortest path joining $\tilde{b}$ and $b$ in $\mathcal{D}$.
	\end{remark}

	\subsection{The proofs of propositions and theorems  given in Sect.~\ref{sect:correctness}}
	%\subsection{Details in then main algorithm}
	\label{app.detail-main-alg}

	%It is a simple matter to obtain the following lemma.
	Lemma~\ref{lem:cut_shortest_path} is presented without a proof (since the proof is similar to that of Proposition 2.19,~\cite{An2017}).
	
	\begin{lemma}
		\label{lem:cut_shortest_path}
		Suppose that $\gamma$ is the shortest path joining two given points $p$ and $q$ in a polygon $\mathcal{A}$ such that $\gamma$ intersects a cutting segment $[u, v]$ at a point, say $a$. If $a \in ]u, v[$, then $\angle \left( \gamma(p,a), \gamma(a,q) \right) = \pi$. %where $\gamma (p,a)$ ($\gamma(a,q)$, respectively) is the sub-path of $\gamma$ from $p$ to $a$ (from  $a$ to $q$, respectively).
	\end{lemma}
	%\begin{proof}
	%	Assume that $\angle \left( \gamma(p,a), \gamma(a,q) \right) \neq \pi$. Then there exists a small-enough  ball $B$ centered at $a$ such that $B \subset \mathcal{D}$ and the intersection of the boundary of $B$ and  $\gamma$ is not empty. Let $s$ and $s'$ be these points of intersection, where $s \in \gamma(p,a)$ and $s'\in \gamma(a,q)$ (see Fig~\ref{fig:cut-polygon-di}(ii)). Let $\gamma_1$, $\gamma_2$ and $\gamma_3$ be the sub-paths of $\gamma$ joining $p$ and $s$; joining $s$ and $s'$; joining $s'$ and $q$, respectively. 
	%	As $\angle \left( \gamma(p,a), \gamma(a,q) \right) \neq \pi$, we obtain $l(\gamma_2) > l([s,s'])$. Therefore we can replace $\gamma$ by  $\gamma_1 \cup [s,s'] \cup \gamma_3$, where the length of latter is strictly less than that of the former, which contradicts the fact that $\gamma$ is the shortest pathục joining $p$ and $q$. Hence $\angle  \left( \gamma(p,a), \gamma(a,q) \right) = \pi$.
	%\end{proof}
	%\begin{figure}[htp]
	%	\centering
	%	\includegraphics[width=0.5\linewidth]{reflex-properties}
	%	\caption{Illustration of the proof the lemma}
	%	\label{fig:reflex-properties}
	%\end{figure}
	
	\medskip
	
\begin{figure}[htp]
		\centering
		\begin{subfigure}[b]{0.28\textwidth}
			\centering
			\includegraphics[scale=0.8]{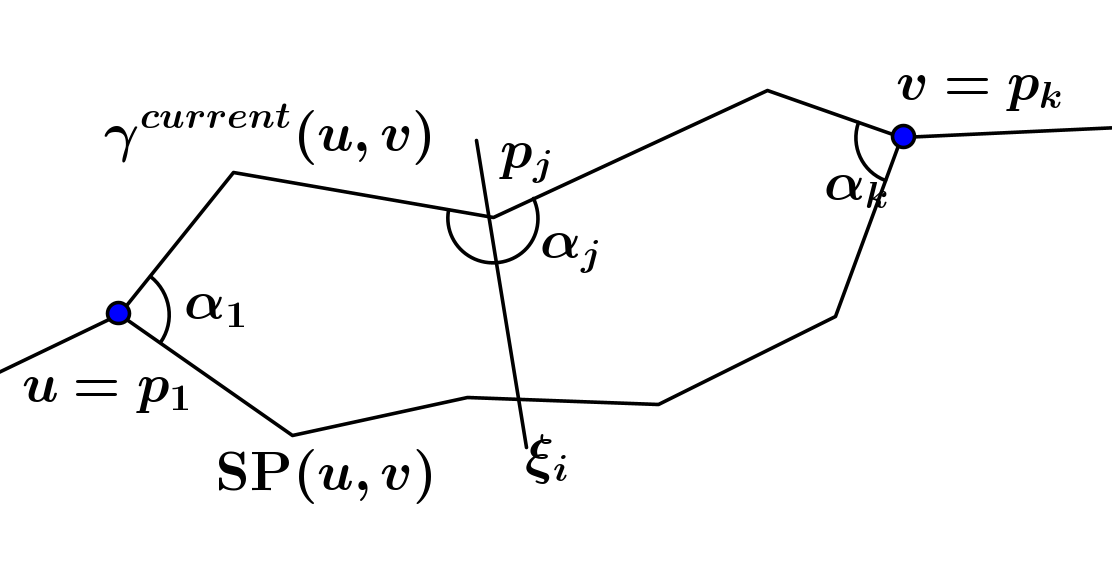}
			\caption*{(i)}
		\end{subfigure}
		\begin{subfigure}[b]{0.7\textwidth}
			%	\centering
			\raggedleft
			\includegraphics[scale=0.7]{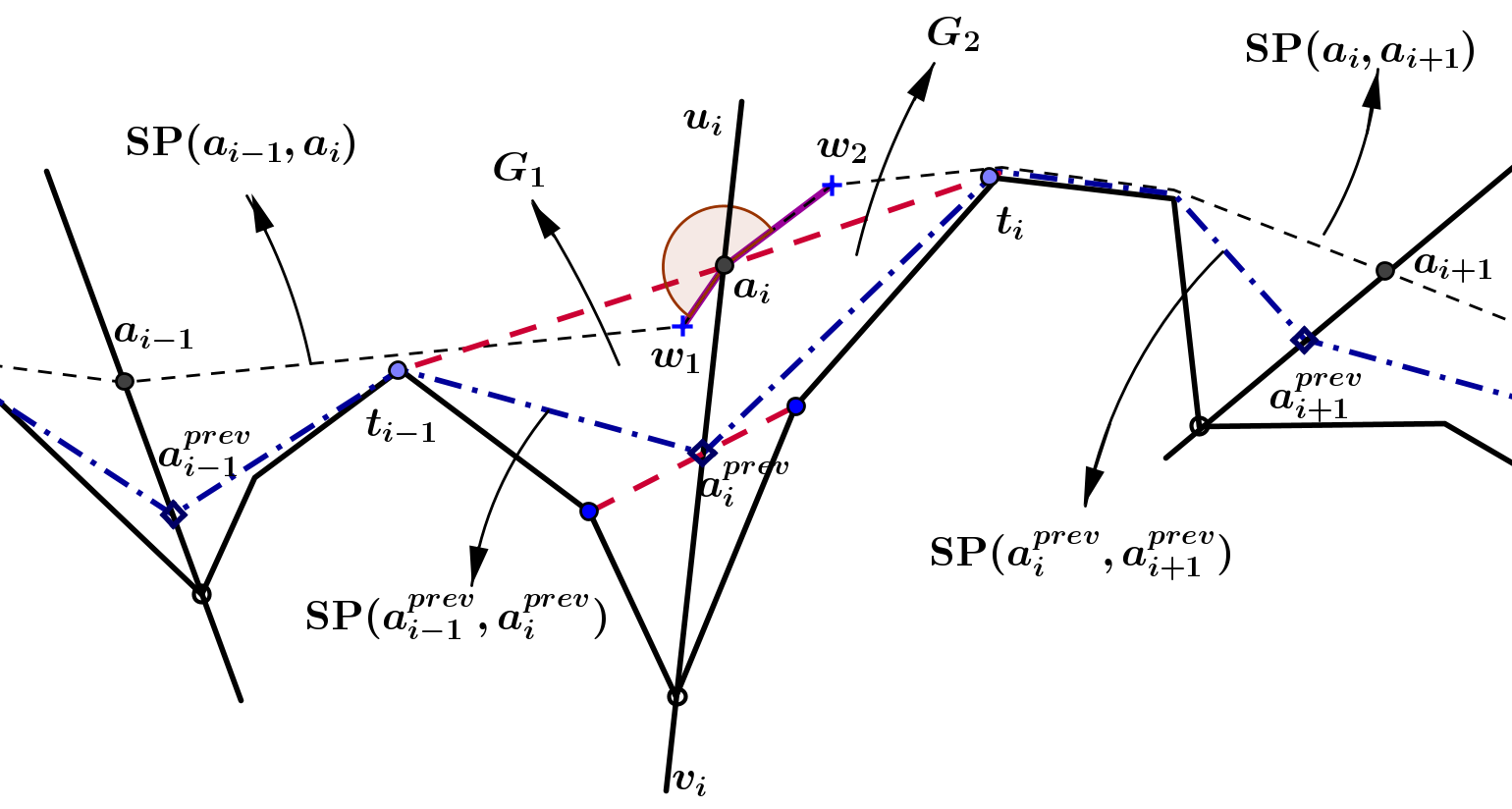}
			\caption*{(ii)}
		\end{subfigure}
		\caption{(i): Illustration of the proof of Proposition~\ref{prop:collinear-condition} and (ii): illustration of the proof of Proposition~\ref{prop:order-gamma-j} in the case $a_i^{prev} \neq a_i$}
		\label{fig:cut-polygon-di}
	\end{figure}

	\noindent \textbf{The proof of Proposition~\ref{prop:collinear-condition}:}
	\begin{proof}
		Let $\gamma^{current} =\cup_{i=0}^{N} \text{SP}(a_i,a_{i+1})$ where $\{a_i\}_{i=1}^N$ is a set of shooting points and  $a_0 = \tilde{b}, a_{N+1} = b$. 	
		Assume the contrary that $\gamma^{current}$ and $\gamma^*$ are distinct. There are two distinct points $u, v$ which are common points to both $\gamma^{current}$ and $\gamma^*$ such that $\gamma^{current} (u,v) \cap \text{SP}(u,v) =\{u,v\}$, where $\gamma^{current}(u,v)$ is the sub-path of $\gamma^{current}$ from $u$ to $v$. Note that $\text{SP}(u,v)$ is also a sub-path of $\gamma^*$. Two points $u$ and $v$ always exist and may be identical to $\tilde{b}$ and $b$, respectively. Then $\gamma^{current}(u,v)$ and $\text{SP}(u,v)$  form a simple polygon with $m$ vertices  $p_1=u, p_2, \ldots, p_k=v, p_{k+1} \ldots, p_m$ (see Fig.~\ref{fig:cut-polygon-di}(i)). Let $\alpha_j$ be the internal angle at $p_j$, with $j=1,2,\ldots, m$. By
		the properties of  shortest paths,
		% Remark~\ref{rem:shortest_path}(a), 
		if $p_j \in \text{SP}(u,v)$ then vertices of  $\text{SP}(u,v)$ except for $u$
		and $v$ are reflex vertices of  $\text{SP}(u,v)$, i.e., $\alpha_j > \pi$.
		If $p_j \in \gamma^{current}(u,v)$, then $p_j$ is either a shooting point or  a reflex vertex of  $\text{SP}(a_i, a_{i+1})$ with some index $i \in \{1,2,\ldots, N\}$. If  $p_j$ is a shooting point, then $\alpha_j \ge \pi$ since Collinear Condition (B) holds at all shooting points. The equality does not occur, because $p_j$ is a vertex of $\gamma^{current}$, and thus $\alpha_j > \pi$.
		% Because $p_j$ is a vertex of $\gamma$, according  to Lemma~\ref{lem:cut_shortest_path},  $\alpha_j > \pi$.
		If $p_j$ is a reflex vertex of  $\text{SP}(a_i, a_{i+1})$ with some index $i \in \{1,2,\ldots, N\}$, we also get $\alpha_j > \pi$.
		%, since proper vertices of the shortest path $\text{SP}(a_i, a_{i+1})$ except for $a_i, a_{i+1}$ are reflex.
		Hence $\alpha_j > \pi$ for $j \in \{1,2,\ldots, m\} \setminus \{1,k\}$.  Since the
		sum of the measures of the internal angles of the simple polygon with $m$ vertices is
		$(m - 2)\pi$, we have 
		$(m-2) \pi =\sum_{j=1}^{m} \alpha_i> \alpha_1 + (m-2) \pi + \alpha_k >(m-2) \pi.$
		The contradiction deduces $\gamma^{current}  = \gamma^*$.
	\end{proof}

	\medskip
	
	\noindent \textbf{The proof of Proposition~\ref{prop:two-stop-condition}:}
	\begin{proof}
		$(\Rightarrow)$ Assuming that Collinear Condition (B) holds  at $a_i$, to show that $\text{SP}(t_{i-1}, t_i) \cap [u_i, v_i] = a_i$, we just need to prove that %
		$\text{SP}(t_{i-1}, t_i)  = \text{SP}(t_{i-1}, a_i) \cup \text{SP}(a_i, t_i)$
		% $\text{SP}(t_{i-1}, t_i)=\gamma^{current}(t_{i-1}, t_i).$ 
		Assume the contrary, we use the same way of the proof of Proposition~\ref{prop:collinear-condition} and the property that the sum of the internal angles in a simple polygon with $m$ vertices is $(m-2)\pi$ to obtain a contradiction. 
		
		$(\Leftarrow)$ Assume that $\text{SP}(t_{i-1}, t_i) \cap [u_i, v_i] = a_i$. 
		%Let $[u,v]$ be the cutting segment corresponding to these shooting points. 
		Therefore $a_i$ belongs to $\text{SP}(t_{i-1}, t_i)$.
		If $a_i = v_i$, then $\angle \left( \text{SP}(t_{i-1}, a_i^{next}), \text{SP}(a_i^{next},t_{i}) \right) \ge \pi$ by the property of shortest paths. If $a_i \in ]u_i,v_i[$, then $\angle \left( \text{SP}(t_{i-1}, a_i^{next}), \text{SP}(a_i^{next},t_{i}) \right) = \pi$ due to Lemma~\ref{lem:cut_shortest_path}.
		Thus	 Collinear Condition (B) holds  at $a_i$. 	  
	\end{proof}
	
	%Proposition~\ref{prop:two-stop-condition} follows the following
	%\begin{remark}
	%	\label{rem:new-shooting-point-properties}
	%	In any cases (Collinear Condition (B) holds or not), the new shooting point $a_i^{next}$ obtained by Procedure \uppercase{Collinear\_Update}($\gamma^{current}, \gamma^{next}, flag$) is  the point of intersection of $\text{SP}(t_{i-1}, t_i)$ and $\xi_i$, for $i=1,2, \ldots,N$.
	%\end{remark}

	%%%% -remove them by replacing the proof in AN HoaiHai 2017%%%%%%%%%%
	%	To prove  Theorem~\ref{theo:global solution} and Theorem~\ref{theo:complexity} we need some results  		 which are not presented in this section. For more detail, see Proposition~\ref{prop:sup-paths-MMS}, Proposition~\ref{prop:order-gamma-j} and Corollary~\ref{coro:inner-angle-shooting-point} in  Appendix (Sect.~\ref{app.detail-main-alg}).
	
	To prove  Theorems~\ref{theo:global solution}, and~\ref{theo:update_path} we need  Proposition~\ref{prop:sup-paths-MMS}, Proposition~\ref{prop:order-gamma-j} and Corollary~\ref{coro:inner-angle-shooting-point}.
	Next, assuming that a partition of $\mathcal{D}$ is given
	by~(\ref{eq:decomposition}), we have 
	\begin{proposition}
		\label{prop:sup-paths-MMS}
		We have	$\gamma^* \subset \cup_{i=0}^{N} \mathcal{D}_i$, where $\mathcal{D}_i$ is bounded by  $B_1,
		B_2$ and the cut  segments $\xi_i$ and $\xi_{i+1}$.  Moreover there exists a set of  points $\{a_i^* \in \xi_i, i=0,1, \ldots, N+1\}$, where $a_0^*=\tilde{b}, a_{N+1}^*=b$, such that $\gamma^*= \cup_{i=0}^{N} {\rm SP}(a_i^*, a_{i+1}^*)$, in which ${\rm SP}(a_i^*, a_{i+1}^*)$ is the shortest path joining $a_i^*$ and $a_{i+1}^*$ in $\mathcal{D}_i$, for $i=0,1, \ldots,N (N \ge 1$).
	\end{proposition}
	\begin{proof}
		Because of (\ref{eq:decomposition}), $\xi_i$ divides $\mathcal{D}$ into two parts, one  contains $\tilde{b}$ and the other  contains $b$. Hence  the intersection of $\gamma^*$ and $\xi_i$ is not empty. As $\gamma^*$ is the shortest path joining $\tilde{b}$ and $b$ in $\mathcal{D}$, the  intersection is only one point. 
		%Indeed, if $\gamma^*$ intersects with $\xi_i$ at more than two distinct points, we can replace $\gamma^*$ by a  path passing through $\xi_i$ at only one of those points, which is shorter, a contradiction. 
		%	Let $a_i^*$ be  the point of intersection of $\gamma^*$ and $\xi_i$, for $i=1,2, \ldots,N$.
		Let $a_i^*=\gamma^* \cap \xi_i$, for $i=1,2, \ldots,N$, where $a_0^*=\tilde{b}, a_{N+1}^*=b$. 
		% As every sub-path joining two arbitrary points of the shortest path is the shortest path,  
		Then $\gamma^*= \cup_{i=0}^{N} \text{SP}(a_i^*, a_{i+1}^*)$, where $\text{SP}(a_i^*, a_{i+1}^*)$ is the shortest path joining $a_i^*$ and $a_{i+1}^*$ in $\mathcal{D}$, for $i=0,1, \ldots,N$.
		Since the construction of $\mathcal{D}_i$, we deduce that $\text{SP}(a_i^*, a_{i+1}^*) \subset\mathcal{D}_i$. 
	\end{proof}  
	
	%\subsubsection*{c) The factor (f2)..}
	The converse of Proposition~\ref{prop:sup-paths-MMS} is not true for the general case. 
	%	It means that if $\gamma$ is a path joining $\tilde{b}$ and $b$ and intersects with $\xi_i$ at $a_i$ such that $\gamma (a_i,a_{i+1})$ is the shortest path joining $a_i$ and $a_{i+1}$ in $\mathcal{D}_i$, for all $i=0,1,\ldots,N$, i.e, $\gamma =\cup_{i=0}^{N} \text{SP}(a_i, a_{i+1})$, then it is not concluded that $\gamma$ is $\gamma^*$. 
	It means that if there is a set $\{a_i \in \xi_i\}_{i=0}^N$,
	$\gamma^{current} =\cup_{i=0}^{N} \text{SP}(a_i, a_{i+1})$, then it is not concluded that $\gamma ^{current}$ is $\gamma^*$. 
	However, if $a_i$ satisfies Collinear Condition (B), for all $i=1,2,\ldots,N$, 
	%	$i=\overline{1,N}$,
 we have $\gamma^{current} = \gamma^*$	by  Proposition~\ref{prop:collinear-condition}. 
	%$\gamma$ is the shortest path joini,g $\tilde{b}$ and $b$ in $\mathcal{D}$.

	%\subsection{The proof of some results}
	%For  $j=0,1, \ldots$, let $\gamma^j$ be the obtained path at $j^{th}$- step of the algorithm,  where $\gamma^j = \cup_{i=0}^N \text{SP}(a_i^j,a_{i+1}^j)$, $a_i^j \in [v_i,u_i[$ is the shooting point at $j^{th}$- step.

	\medskip
	
	%	According to the properties of the shortest paths and Remark~\ref{rem:shortest_path}(a),
	% and Remark~\ref{rem: path-crosses}),  	we obtain the following 
	\begin{remark}
		\label{rem:two-non-cross-SP}
		%	Let $\mathcal{D}_i$ be the region bounded by $B_1, B_2$, cutting segments $\xi_i$ and $\xi_{i+1}$. 
		Suppose that $a_i \in \xi_i, a_{i+1} \in \xi_{i+1}$ are shooting points in some iteration step, and $a_i^{next} \in \xi_i, a_{i+1}^{next} \in \xi_{i+1}$ are shooting points at  the next iteration step of the proposed algorithm such that $a_i^{next} \in [a_i, u_i]$ and $a_{i+1}^{next} \in [a_{i+1}, u_{i+1}]$, see Figs.~\ref{fig:MMS}(iii) and (iv).
		%Let $\tau_1$ ($\tau_2$, resp.) be  the shortest path joining $a_i$ and $a_{i+1}$ ($a_{i}^{next}$ and $a_{i+1}^{next}$, resp.) in $\mathcal{D}_i$.
		%, see Fig~\ref{fig:cut-polygon-di}(ii). 
		Then we have
		\begin{itemize}
			% 			\item[(a)]  $\tau_1$ and $\tau_2$ are convex  with convexity facing towards $B_1$.
			% 			\item[(b)]   $\tau_2$ is entirely contained in the polygon bounded by $[a_i^{next}, a_i]$, $\tau_1$, $[a_{i+1}, a_{i+1}^{next}]$ and $B_1$. Thus we say that  $\tau_2$  is above $\tau_1$ when we view from $B_1$.
			\item[(a)]  ${\rm SP}(a_i,a_{i+1})$ and ${\rm SP}(a_i^{next},a_{i+1}^{next})$ are convex  with convexity facing towards $B_1$.
			\item[(b)]   ${\rm SP}(a_i^{next},a_{i+1}^{next})$ is entirely contained in the polygon bounded by $[u_i, a_i]$, ${\rm SP}(a_i,a_{i+1})$, \break $[a_{i+1}, u_{i+1}]$, and $B_1$. Thus we say that  ${\rm SP}(a_i^{next},a_{i+1}^{next})$  is above ${\rm SP}(a_i,a_{i+1})$ when we view from $B_1$.
		\end{itemize}
	\end{remark}

	% 	\begin{figure}[htp]
	% 	\centering
	% 	\begin{subfigure}[b]{0.48\textwidth}
	% 		\centering
	% 		\includegraphics[scale=0.68]{proof-prop1}
	% 		\caption*{(i)}
	% 	\end{subfigure}
	% 	\begin{subfigure}[b]{0.48\textwidth}
	% 		\centering
	% 		\includegraphics[scale=0.45]{cut-polygon-Di}
	% 		\caption*{(ii)}
	% 	\end{subfigure}
	% 	\caption{(i)- the illustration of the proof of Proposition~\ref{prop:collinear-condition} and (ii)- the illustration of Remark~\ref{rem:two-non-cross-SP}}
	% 	\label{fig:cut-polygon-di}
	% \end{figure}

	%	\begin{figure}[htp]
	%	\centering
	%		\includegraphics[scale=0.65]{cut-polygon-Di}
	%	\caption{Illustration of Remark~\ref{rem:two-non-cross-SP}}
	%	\label{fig:cut-polygon-di}
	%\end{figure}

	\begin{lemma}
		\label{lem:gamma_j}
		Recall that the upper angle created by ${\rm SP}(a_{i-1},a_i)$ and ${\rm SP}(a_i,a_{i+1})$ w.r.t. $[u_i, v_i]$ is denoted by $\angle \left( {\rm SP}(a_{i-1},a_i), {\rm SP}(a_i,a_{i+1}) \right)$. For $i=1,2,\ldots,N$, we suppose that $a_i \in ]u_i,v_i[$. Let $\hat{a_i} = \text{SP}(t_{i-1}, t_i) \cap [u_i, v_i]$, where $t_{i-1}$ and $t_i$ are temporary points w.r.t. $a_{i-1}$ and $a_i$.
		Then we have
		\begin{itemize}
			\item[(a)]  $\angle \left( {\rm SP}(a_{i-1},a_i), {\rm SP}(a_i,a_{i+1}) \right) = \pi \Leftrightarrow \hat{a_i} = a_i$;
			\item[(b)]  $\angle \left( {\rm SP}(a_{i-1},a_i), {\rm SP}(a_i,a_{i+1}) \right) < \pi \Leftrightarrow \hat{a_i} \in ]a_i, u_i[$;
			\item[(c)]  $\angle \left( {\rm SP}(a_{i-1},a_i), {\rm SP}(a_i,a_{i+1}) \right) > \pi \Leftrightarrow \hat{a_i} \in [v_i, a_i[$.
		\end{itemize}
	\end{lemma}
	\begin{proof}
		When Collinear Condition (B) does not hold, the role of $\hat{a_i}$ and $a_i^{next}$ are the same.
		By Proposition~\ref{prop:two-stop-condition}, Collinear Condition (B) in Sect.~\ref{sect:collinear_condition}, and the fact that $a_i \neq v_i$, we get (a). % We will use a proof by contradiction for (b) and (c) as follows
		
		b) $(\Rightarrow)$	  Suppose that $\angle \left( \text{SP}(a_{i-1},a_i), \text{SP}(a_i,a_{i+1}) \right)< \pi$ 		and  $\hat{a_i} \notin ]a_i, u_i]$. 
		Because $\text{SP}(t_{i-1}, t_i)$ never touches $B_1$ of $\mathcal{D}$, similarly to Remark~\ref{remark:boundaryB_2}, then $\hat{a_i} \neq u_i$.
		%, where $t_{i-1}$ and $t_i$ are temporary points w.r.t. $a_{i-1}$ and $a_i$. 
		Therefore $\hat{a_i} \in [v_i,a_i]$. 
		%
		%	Since $\angle \left( \text{SP}(a_{i-1},a_i), \text{SP}(a_i,a_{i+1}) \right)< \pi$, 
		Clearly, the collinear condition does not hold at $a_i$.		Therefore $a_i \neq \hat{a_i}$,  $\text{SP}(t_{i-1}, t_i)$ and $\text{SP}(t_{i-1},a_i) \cup \text{SP}(a_i,t_{i})$ are thus distinct. Then there are two distinct points $u,v$ which
		are common points to both $\text{SP}(t_{i-1}, t_i)$ and $\text{SP}(t_{i-1},a_i) \cup \text{SP}(a_i,t_{i})$ such that 
		$\text{SP}(u,v) \cap \left(\text{SP}(u,a_i) \cup \text{SP}(a_i,v) \right) = \{u, v\}.$
		Note that $\text{SP}(u,v)$ is also a sub-path of $\text{SP}(t_{i-1}, t_i)$.
		Then $\text{SP}(u,v)$ and $\text{SP}(u,a_i) \cup \text{SP}(a_i,v)$  form a simple polygon in $\mathcal{D}$ with $m$ vertices $p_1=u, p_2, \ldots, p_l=a_i, p_{l+1} \ldots, p_k=v, p_{k+1}, \ldots,  p_m$. Let $\alpha_j$ be the internal angle at $p_j$, with $j=1,2,\ldots, m$, see Fig~\ref{fig:cut-polygon-di}(i). According to the property of shortest paths, $\alpha_j >\pi$ for $j \in \{1,2,\ldots, m\} \setminus \{1,k,l\}$.
		If $j=l$,  due to $\alpha_j = 2\pi - \angle \left( \text{SP}(a_{i-1},a_i), \text{SP}(a_i,a_{i+1}) \right)$, we have $\alpha_j> \pi$. Then $\alpha_j >\pi$ for $j \in \{1,2,\ldots, m\} \setminus \{1,k\}$.
		%
		%Due to Remark~\ref{rem:shortest_path}(a), 
		%vertices of  shortest paths except source and destination points are reflex vertices of  the paths, i.e., 
		% $\alpha_s \ge \pi$, for $s \in \{1,2, \ldots, m\} \setminus \{1,k,l\}$.
		%
		%	In the case $p_j$ is a vertex of  $\text{SP}(a_i, a_{i+1})$ with some index $i \in \{1,2,\ldots, N\}$, we also get $\alpha_j > \pi$ since vertices of the shortest path $\text{SP}(a_i, a_{i+1})$ except for $a_i, a_{i+1}$ are reflex. 		Hence $\alpha_i > \pi$ with $i \in \{1,2,\ldots, m\} \setminus \{1,j\}$.
		%	We have $\alpha_k = 2\pi - \angle \left( \text{SP}(a_{i-1},a_i), \text{SP}(a_i,a_{i+1}) \right)$ and thus   $\alpha_k> \pi$.
		%		Because the 		sum of the measures of the internal angles of the simple polygon with $m$ vertices is
		%		$(m - 2)\pi$, we have the following
		%		$$(m-2) \pi =\sum_{s=1}^{m} \alpha_s \ge  (m-3) \pi + \alpha^-+ \alpha_k+ \alpha_l\ge  (m-3) \pi + \alpha^-+ \pi + \alpha_l  >(m-2) \pi.$$
		It is similar to the proof of Proposition~\ref{prop:collinear-condition}, a  contradiction is obtained. Thus $\hat{a_i} \in ]a_i, u_i[$.

		c) $(\Rightarrow)$			Suppose that $\angle \left( \text{SP}(a_{i-1},a_i), \text{SP}(a_i,a_{i+1}) \right) > \pi$ and  $\hat{a_i} \notin ]v_i, a_i[$. A contradiction can be deduced in much the same way as the above proof. 
		Therefore $\hat{a_i} \in ]v_i, a_i[$.
		
		%Thus the  necessary conditions of (b) and (c)  are obtained. For the case of sufficient condition of (b), ((c) can be considered similarly), assume that $\hat{a_i} \in ]a_i, u_i[$. 
		b) c) $(\Leftarrow)$ For the sufficient condition of (b), ((c) can be considered similarly), assume that $\hat{a_i} \in ]a_i, u_i[$. 		Then 
		$\angle \left( \text{SP}(a_{i-1},a_i), \text{SP}(a_i,a_{i+1}) \right) < \pi$, for if not, $\angle \left( \text{SP}(a_{i-1},a_i), \text{SP}(a_i,a_{i+1}) \right) \ge \pi$ deduces that $\hat{a_i} \in [v_i, a_i]$ by the  necessary conditions of (a) and (c) that are proved above, which is impossible. Hence the proof is complete. 
		%		\item[(b)] 
		%		Suppose that  $\angle \left( \text{SP}(a_{i-1},a_i), \text{SP}(a_i,a_{i+1}) \right) = \pi$ and $\hat{a_i} = a_i$. The contradiction may be deduced in much the same way as the proof of (a) and note that $\alpha_k = \angle \left( \text{SP}(a_{i-1},a_i), \text{SP}(a_i,a_{i+1}) \right)= \pi$.
		%		\item[(c)] Suppose that $\angle \left( \text{SP}(a_{i-1},a_i), \text{SP}(a_i,a_{i+1}) \right) > \pi$ and $\hat{a_i} \notin [v_i, a_i[$. Analysis similar to that in the proof of (a) and note that $\alpha_k = \angle \left( \text{SP}(a_{i-1},a_i), \text{SP}(a_i,a_{i+1}) \right) > \pi$, we also obtain a contradiction.
	\end{proof}
	
	%Because of initial shooting points $a_i^0 =v_i$, for $i=1,2,\ldots,N$, we have the following 
	%\begin{remark}
	%	\label{rem:gamma_0}
	%	By taking initial shooting points $a_i^0 =v_i$, for $i=1,2,\ldots,N$ and according to Procedure \uppercase{Collinear\_Update}($\gamma^{current}, \gamma^{next}, flag$), we obtain $a_i^1 \in [a_i^0, u_i[$, where $a_i^0$ and $a_i^1$  are the shooting points at initial step and $1^{th}$-step respectively, for $i=1,2,\ldots,N$.
	%\end{remark}

	\begin{proposition}
		\label{prop:order-gamma-j}
		
		For all $j$, we have $\gamma^{j+1}$ is above $\gamma^{j}$ when we view from $B_1$, i.e., $a_i^{j+1} \in [a_i^j, u_i[$, for $i=1,2, \ldots, N$, where $\gamma^{j}$ is the path obtained in $j^{th}$-iteration step of the algorithm.

	\end{proposition}
	
	\begin{proof}
		We give a proof by induction on $j$.
		%	It is clear that 
		The statement holds for $j=1$ due to taking initial shooting points $a_i^0 =v_i$, for $i=1,2,\ldots,N$.
		Let $\{a_i^{prev}\}_{i=1}^N, \{a_i\}_{i=1}^N$ and $\{a_i^{next}\}_{i=1}^N$ be the sets of shooting points corresponding to $k-1^{th}$, $k^{th}$ and $k+1^{th}$-iteration steps of the algorithm, respectively ($k \ge 1$).
		Assuming that,  for $i=1,2, \ldots, N$,  $a_i \in [a_i^{prev}, u_i[$, we next prove that $a_i^{next} \in [a_i, u_i[$.
		
		If  $\angle \left( \text{SP}(a_{i-1},a_i), \text{SP}(a_i,a_{i+1}) \right) \le \pi$ due to  Lemma~\ref{lem:gamma_j}(a) and (b), we  obtain  $a_i^{next} \in [a_i, u_i[$. Thus we just need to prove that $a_i^{next} \in [a_i, u_i[$ in the case $\angle \left( \text{SP}(a_{i-1},a_i), \text{SP}(a_i,a_{i+1}) \right) > \pi$.
		%, for $i=1,2, \ldots, N$.
		If $a_i = v_i$, then Collinear Condition (B) holds at  $a_i$ then $a_i^{next} = a_i \in [a_i, u_i[$.
		If $a_i \neq v_i$,  suppose contrary to our claim, there is an index $i$ such that $\angle \left( \text{SP}(a_{i-1},a_i), \text{SP}(a_i,a_{i+1}) \right) > \pi$ and $a_i^{next} \notin [a_i, u_i[$, i.e., 
		\begin{align}
			\label{eq:7}
			a_i^{next} \in [v_i, a_i[ 
		\end{align}

		Let $G_1$  be the closed region bounded by $[a_{i}^{prev},a_i]$, $\text{SP}(a_i,t_{i-1})$ and $\text{SP}(t_{i-1},a_{i}^{prev})$. Let $G_2$  be the closed region bounded by $[a_{i}^{prev},a_i]$, $\text{SP}(a_i, t_{i})$ and $\text{SP}(t_{i},a_{i}^{prev})$, where $t_{i-1}$ and $t_{i}$ are   temporary points w.r.t. $a_{i-1}^{prev}$ and  $a_i^{prev}$ (see Fig.~\ref{fig:cut-polygon-di}(ii)). Since $a_i \in ]u_i,v_i[$, $\angle \left( \text{SP}(t_{i-1},a_i), \text{SP}(a_i,t_{i}) \right) = \pi$ by Lemma~\ref{lem:cut_shortest_path}.
		%	If $a_i = a_i^{prev}$, $G_1$ and $G_2$ are polylines. If $a_i \neq a_i^{prev}$, by properties of shortest paths, $G_1$ and $G_2$ have the form of a funnel. This concept is introduced explicitly in~\cite{An2017, Lee1984}.
		%
		Next, let $[a_i,w_1]$ and $[a_i, w_2]$ be segments having
		the common endpoint $a_i$ of $\text{SP}(a_{i-1},a_i)$ and $\text{SP}(a_i,a_{i+1})$, respectively.
		As  $\angle \left( \text{SP}(t_{i-1},a_i), \text{SP}(a_i,t_{i}) \right) = \pi$ and  $\angle \left( \text{SP}(a_{i-1},a_i), \text{SP}(a_i,a_{i+1}) \right) > \pi$,  there are at least one of two following cases that will happen: 
		\begin{align}
			\label{eq:two-cases-cross}
			&	]a_i, w_1]\text{ is below entirely  } \text{SP}(t_{i-1},a_i) \\ \notag
			&	\text{    or } ]a_i, w_2] \text{  is below entirely } \text{SP}(a_i, t_{i}) \text{   when we view from } u_i. \hspace{4cm}
		\end{align} 
		%
		% 		\begin{figure}[htp]
		% 			\centering
		% 			\includegraphics[width=0.7\linewidth]{proof_main_prop}
		% 			\caption{Illustration of the proof of the proposition in the case $a_i^{prev} \neq a_i$}
		% 			\label{fig:proofmainprop}
		% 		\end{figure}
		
		%W.l.o.g, assume that $]a_i, w_2]$ is below $\text{SP}(a_i, t_{i})$. 
		%	Note that $\text{SP}(a_i, t_{i})$ is also a shortest path.  
		Then $\text{SP}(a_i, t_{i})$, $\text{SP}(a_i,a_{i+1})$ and $\text{SP}(a_i^{prev},t_i)$ do not cross each other, due to the properties of shortest paths.
		According  to  Remark~\ref{rem:two-non-cross-SP}(a),  for any temporary point $t_i \in \text{SP}(a_i^{prev},a_{i+1}^{prev})$, 
		there is a polygon bounded by $[a_i^{prev}, a_i]$,  $[a_{i+1}, a_{i+1}^{prev}]$, $\text{SP}(a_i,a_{i+1})$, and $\text{SP}(a_i^{prev},a_{i+1}^{prev})$ which contains entirely $\text{SP}(a_i, t_{i})$. Thus $\text{SP}(a_i,a_{i+1})$ does not intersect with the interior of $G_2$. Similarly, $\text{SP}(a_{i-1}, a_i)$ does not also intersect with the interior of $G_1$. These things contradict~(\ref{eq:two-cases-cross}), which  completes the proof.
	\end{proof}

	%	The value of the upper angle at a shooting point is established by our next corollary, which is clear from Proposition~\ref{prop:order-gamma-j} and Lemma~\ref{lem:gamma_j}.
	
	\begin{corollary}
		\label{coro:inner-angle-shooting-point}
		Let $\gamma^{current}$ be the  path obtained in some iteration step of the algorithm, where $\gamma^{current} = \cup_{i=0}^N {\rm SP}(a_i,a_{i+1})$.
		If $a_i \in ]v_i, u_i[$, then  we have $\angle \left( {\rm SP}(a_{i-1},a_i), {\rm SP}(a_i,a_{i+1}) \right) \le \pi$, for $i=1,2, \ldots, N$.
	\end{corollary}
	
	\begin{lemma}[Proposition 5.1, \cite{AnHaiHoai2013}]
		\label{lem:AnHaiHoai}
		If $a , b , c,$ and $d$ are points in a simple polygon $\mathcal{P}$ such that $[ a , d ], [ b , c ] \subset \mathcal{P}$ , then
		$d_{\mathcal{H}} ({\rm SP}(a,b), {\rm SP}(c,d)) \le \Vert a-d \Vert + \Vert b-c \Vert \le 2 \max \{ \Vert a-d \Vert, \Vert b-c \Vert \}.$
	\end{lemma}
	
	\begin{lemma}[Lemma 3.1, \cite{HaiAn2011}]
		\label{lem:HaiAn2011}
		If $\{a_n\}$ and $\{b_n\}$ are sequences of points in a simple polygon, $a_n \rightarrow a$ and
		$ b_n \rightarrow b$, then
		$d_{\mathcal{H}} ({\rm SP}(a_n,b_n), {\rm SP}(a, b)) \rightarrow 0$ as $n \rightarrow +\infty.$
	\end{lemma}
	
	% {\color{red}
	% 	\begin{lemma}[Theorem 3.1 in \cite{HaiAn2011}]
	% 	\label{theo:HaiAn2011}
	% Suppose
	% that $\mathcal{A}$ is a simple polygon.
	% If $A$ is a nonempty compact subset of $\mathcal{A}$ and $d_{\mathcal{H}} ({\rm SP}(a_n,b_n), A) \rightarrow 0$ as $n \rightarrow +\infty,$
	% then $A$ is a shortest path.
	% \end{lemma}
	% }
	% 	
	\medskip
	\noindent \textbf{The proof of Theorem~\ref{theo:global solution}:}
	\begin{proof}
		Let $\gamma^j = \cup_{i=0}^N \text{SP}(a_i^j,a_{i+1}^j)$, where $a_i^j \in [u_i, v_i]$ for $i=1,2, \ldots,N$ and $a_0^j = \tilde{b}, a_N^j =b$. 
		If the algorithm stops after a finite number of iterations, the path obtained  is the shortest. We thus just need to consider the case that $\{\gamma^j\}_{j\in \mathbb{N}}$ is infinite, according to Proposition~\ref{prop:collinear-condition}.

		For each $i=1,2,\ldots, N$, 
		%since $\{\gamma^j\}_{j\in \mathbb{N}}$ is infinite, w.l.o.g we can assume that $\{a_i^j\}_{j\in \mathbb{N}}$ is infinite.
		%
		since $[u_i,v_i]$ is compact  and $\{a_i^j\}_{j\in \mathbb{N}} \subset [u_i,v_i]$ there exists a sub sequence $\{a_i^{j_k}\}_{k \in \mathbb{N}} \subset \{a_i^j\}_{j\in \mathbb{N}}$ such that $\{a_i^{j_k}\}_{k \in \mathbb{N}}$ converges to a point, say $a_i^*$, in $\mathbb{R}^2$ with $\Vert . \Vert$ as $k \rightarrow \infty$. 
		As $[u_i, v_i]$ is closed, $\{a_i^{j_k}\}_{k \in \mathbb{N}} \subset [u_i,v_i]$, and $\{a_i^{j_k}\}_{k \in \mathbb{N}} \rightarrow a_i^*$, we get $a_i^* \in [u_i, v_i]$. 
		Write $a_0^*=\tilde{b}$ and $a_{N+1}^*=b$. Set $\hat{\gamma} = \cup_{i=0}^N \text{SP}(a_i^*,a_{i+1}^*)$.

		\textbf{i.} \quad \textit{We begin by proving the convergence of whole sequences,i.e., $\{a_i^j\}_{j\in \mathbb{N}} \rightarrow a_i^*$ as $j \rightarrow \infty$, for $i=1,2, \ldots, N$, based on the order of elements of $\{a_i^j\}_{j\in \mathbb{N}}$ shown in Proposition~\ref{prop:order-gamma-j} and the convergence of their subsequence.}
		By the order of elements of $\{a_i^j\}_{j\in \mathbb{N}}$ as shown in Proposition~\ref{prop:order-gamma-j}, for all natural numbers large enough  $j$, there exits $k \in \mathbb{N}$ such that $ a_i^j \in [a_i^{j_k},a_i^{j_{k+1}}]$.
		Furthermore,  we also obtain $\{a_i^{j_k}\}_{k\in \mathbb{N}}$ converges on one side to $a_i^*$.
		%	Since $\{a_i^{j_k}\}_{k \in \mathbb{N}} \rightarrow a_i^*$ as $k \rightarrow \infty$, 
		For all $\delta >0$, there exists $k_0 \in \mathbb{N}$ such that $\Vert a_i^{j_k} - a_i^* \Vert < \delta$, for $k \ge k_0$.
		%
		%	Note that $\{a_i^{j_k}\}_{k \in \mathbb{N}} \subset \{a_i^j\}_{j\in \mathbb{N}}$, $\{a_i^{j_k}\}_{k \in \mathbb{N}} \rightarrow a_i^*$ as $k \rightarrow \infty$ and by Proposition~\ref{prop:order-gamma-j}, we obtain $a_i^{j_k+1} \in [a_i^{j_k }, a_i^*]$ for $j_k \ge 0$ and $\{a_i^{j_k}\} \subset [v_i, a_i^*]$.
		As $a_i^{j} \in [a_i^{j_{k_0}}, a_i^*]$, for $j \ge j_{k_0}$, we get $\Vert a_i^{j} - a_i^* \Vert < \delta$, for $j \ge j_{k_0}$. 
		This clearly forces $\{a_i^{j}\}_{j \in \mathbb{N}} \rightarrow a_i^*$ as $j \rightarrow \infty$. This  is the one-sided convergence, $\gamma^j$ is thus below $\hat{\gamma}$ when we view from the boundary $B_1$.

		\textbf{ii.} \quad \textit{We next indicate that  $d_{\mathcal{H}}(\gamma^j, \hat{\gamma}) \rightarrow 0$ as $j \rightarrow \infty$.}
		According to item \textbf{i},  $a^j_i \rightarrow a_i^*$ as $j \rightarrow \infty$, for $i=1,2,\ldots,N$. By Lemma~\ref{lem:HaiAn2011}, we get $d_{\mathcal{H}}\left( {\rm SP}(a^j_i,a^j_{i+1}),{\rm SP} (a_i^*,a_{i+1}^*) \right) \rightarrow 0$ as $j \rightarrow \infty$.  
		Note that $d_{\mathcal{H}} (A \cup B, C \cup D) \le \max \{ d_{\mathcal{H}} (A,	C), d_{\mathcal{H}} (B,D)\}$, for all closed sets $A,B,C$ and $D$ in $\mathbb{R}^2$, we have
		$d_{\mathcal{H}} (\gamma^j, \hat{\gamma}) \le \max_{0 \le i \le N} \{d_{\mathcal{H}}\left( {\rm SP}(a^j_i,a^j_{i+1}),{\rm SP} (a_i^*,a_{i+1}^*) \right) \}.$
		Therefore $d_{\mathcal{H}}(\gamma^j, \hat{\gamma}) \rightarrow 0$ as $j \rightarrow \infty$.
		%	Then  Remark~\ref{rem:closed-limit}(c) shows that $\lim\limits_{j\rightarrow \infty} \gamma^j = \hat{\gamma}$.

		\textbf{iii.} \quad \textit{Our next claim is that $\hat{\gamma}$ is the shortest path joining $\tilde{b}$ and $b$.} According to Proposition~\ref{prop:collinear-condition}, we need to prove that Collinear Condition (B)  holds at all $a_i^*$, for $i \in \{1,2, \ldots, N\}$.
		Conversely, suppose that  there exists an index $i \in \{1,2, \ldots, N\}$ such that Collinear Condition (B) does not hold at $a_i^*$.
		If $a_i^* = v_i$ then $a_i^j = v_i$ for all $j \ge 0$. Then  Collinear Condition (B) holds at $a_i^j$ for all $j \ge 0$, i.e., 
		$\angle \left( \text{SP}(a_{i-1}^j,a_i^j), \text{SP}(a_i^j,a_{i+1}^j) \right) \ge \pi$, for all $j \ge 0$. Since $\gamma^j$ is  below $\hat{\gamma}$ when we view from the boundary $B_1$ and $\{a_{i-1}^{j}\}_{j \in \mathbb{N}} \rightarrow a_{i-1}^*$,  $\{a_{i+1}^{j}\}_{j \in \mathbb{N}} \rightarrow a_{i+1}^*$ as $j \rightarrow \infty$, 
		we get the upper angle of $\hat{\gamma}$ at $a_i^*$ is not less than $\pi$, which contradicts to Collinear Condition (B) not holding at $a_i^*$.
		Thus $a_i^* \neq v_i$. 
		As Collinear Condition (B) does not hold at $a_i^*$, we conclude $\angle \left( \text{SP}(a_{i-1}^*,a_i^*), \text{SP}(a_i^*,a_{i+1}^*) \right) < \pi$ or $\angle \left( \text{SP}(a_{i-1}^*,a_i^*), \text{SP}(a_i^*,a_{i+1}^*) \right) > \pi.$
		%To obtain a contradiction, we consider two cases

		\textbf{\textit{Case 1:}} $\angle \left( \text{SP}(a_{i-1}^*,a_i^*), \text{SP}(a_i^*,a_{i+1}^*) \right) < \pi$.
		\textit{To obtain a contradiction, we will show that there exists a natural number $j_0$ such that the updating $a_i^{j_0}$ to $a_i^{j_0+1}$ gives  $a_i^{j_0+1} \in (a_i^*, u_i]$.}
		Let $V$ be the set of all vertices of $\mathcal{D}$ which do not belong to $\hat{\gamma}$. Set 
		\begin{align}
			\label{eq:take-epsilon-0}
			\mu = \min \left\{  d (v, \hat{\gamma}): v \in V \right\}.
		\end{align}
		
		Because $V$ is finite, $d (v, \hat{\gamma})>0$, for $v \in V $,  we get $\mu >0$.
		Let $\epsilon_0 = \mu /4$. % be a very small positive number satisfying $\epsilon_0 < \mu$.
		%
		%There is no vertex of $\mathcal{D}$ contained in $N(\hat{\gamma}, \epsilon_0) \setminus \{\hat{\gamma}\}$ according to~(\ref{eq:take-epsilon-0}).
		%
		Since  $a_i^* \neq v_i$, there is  a point in $[v_i, a_i^*]$, say $c_i$, such that $\Vert c_i - a_i^* \Vert = \epsilon_0$ (see Fig.~\ref{fig:prooftheorem}(i)).
		Through $c_i$, we draw a polyline $\mathcal{L}$ whose line segments are parallel to line segments of $\text{SP}(a_{i-1}^*,a_i^*) \cup \text{SP}(a_i^*,a_{i+1}^*)$, respectively (see Fig.~\ref{fig:prooftheorem}(ii)). 
		%Then $\mathcal{L} \subset \mathcal{D}$ and $\mathcal{L}$ does not contain any vertex of $\mathcal{D}$ except vertices belonging to $\hat{\gamma}$.
		Let us denote by $c_{i-1}$ ($c_{i+1}$, resp.) the point of intersection of $\mathcal{L}$ and the line going through $[u_{i-1}, v_{i-1}]$ ($[u_{i+1}, v_{i+1}]$, resp.). 
		%
		% Since $a_i^*$ is contained in the closed ball centered at $c_i$ and radius $\epsilon_0$ which is a very small positive number, w.l.o.g and for simplicity, we can suppose that there is a line going through $a_i^*$ and creating with $\mathcal{L}$ an isosceles triangle, say $\triangle c_idd'$ (see Fig.~\ref{fig:prooftheorem}(i)). 
		%	In the general case, when the polygon created by  a fixed line going through $a_i^*$ and $\mathcal{L}$ is not an isosceles triangle, and the polyline joining $c_i$ and $d$ (or $d'$) is then not a line segment. However we can calculate the lengths of them.
		Since $a_i^*$ is very closed to $c_i$, w.l.o.g we can suppose that there is a line going through $a_i^*$ and creating with $\mathcal{L}$ an triangle, say $\triangle c_idd'$ (see Fig.~\ref{fig:prooftheorem}(i)). If not, we can take
		\begin{align}
			\label{eq:arbitrarily small}
			\epsilon_0 = \dfrac{\mu}{2^{K_0}}, \text{ where } K_0 \text{ is a large enough natural number.}
		\end{align}
		For simplicity, assume that  $\triangle c_idd'$ is an isosceles triangle, (see Fig.~\ref{fig:prooftheorem}(i)).

		%	We next find the lengths of $[a_i,d]$ and $[c_i,d']$. 
		
		\textbf{Claim 1.}\textit{ Let $ \alpha^- =\angle \left( \text{SP}(a_{i-1}^*,a_i^*), [a_i^*,u_i] \right)$ and  $ \alpha^+ =\angle \left( [a_i^*,u_i], \text{SP}(a_i^*,a_{i+1}^*) \right)$. In $\triangle c_idd'$, we obtain}
		\begin{align}
			\label{eq:laws-of-sines}
			\Vert c_i - d \Vert = \epsilon_0 . \dfrac{\sin \left( \frac{1}{2} (\pi -\alpha^- - \alpha^+) \right)}{\sin \left( \frac{1}{2} (\pi -\alpha^- + \alpha^+) \right)}.
		\end{align} 
		\begin{proof}[The proof of Claim 1] Using the Law of Sines in a triangle.
		\end{proof}

		Turning to the proof, for $i = 1,2, \ldots,N$, since $\{a_i^j\}_{j\in \mathbb{N}} \rightarrow a_i^*$ as $j \rightarrow \infty$, there exists $\overline{j} \in \mathbb{N}$ such that   
		\begin{align}
			\label{eq:parameter-M}
			\Vert a_i^j -a_i^* \Vert < \dfrac{\epsilon_0}{M}, \text{ for all } j \ge \overline{j},
		\end{align}
		where $M$ is a given real number which is chosen as in  Claim 2 below and note that $M$ only depends on the segments $[u_{i-1},v_{i-1}], [u_i,v_i]$, and $[u_{i+1},v_{i+1}]$.
		
		\medskip
		\textbf{Claim 2.} 
		\textit{If $[u_{i-1},v_{i-1}]$  is parallel to $[u_i,v_i]$, then we set $m_1 = 1$. Otherwise let the point of intersection of two lines going through  $[u_{i-1},v_{i-1}]$ and $[u_i,v_i]$  be $s_1$ and we set $m_1 = \dfrac{\Vert s_1 - a_i^* \Vert}{\Vert s_1 - a_{i-1}^* \Vert}$.
			Similarly, if $[u_{i+1},v_{i+1}]$  is parallel to $[u_i,v_i]$, then we set $m_2 = 1$. Otherwise let the point of intersection of two lines going through $[u_{i+1},v_{i+1}]$ and $[u_i,v_i]$  be $s_2$ and we set $m_2 = \dfrac{\Vert s_2 - a_i^* \Vert}{\Vert s_2 - a_{i+1}^* \Vert}$.
			Take $M = \max \{1, m_1, m_2 \}$.
			Then  ${\rm SP}(a_{i-1}^{j},a_i^{j})$ and ${\rm SP}(a_i^{j},a_{i+1}^{j})$ are contained in the region bounded by $\mathcal{L}, [c_{i+1}, a_{i+1}^*], {\rm SP}(a_{i-1}^*,a_i^*) \cup {\rm SP}(a_i^*,a_{i+1}^*)$, and $[a_{i-1}^*, c_{i-1}]$, for all $j \ge \overline{j}$, where $\overline{j}$ is given by (\ref{eq:parameter-M}). }
		
		\begin{proof}[The proof of Claim 2]
			For simplicity, we only prove the case that $[u_{i-1},v_{i-1}]$  is parallel to $[u_i,v_i]$ and two lines going through $[u_{i+1},v_{i+1}]$ and $[u_i,v_i]$ intersect at $s_2$ (see Fig.~\ref{fig:prooftheorem}(ii)).
			One sees immediately that $m_1 =1$ and $\Vert  c_{i-1} - a_{i-1}^*\Vert = \Vert c_{i} - a_{i}^*\Vert = \epsilon_0$. 
			By (\ref{eq:parameter-M}), we get $ \Vert  a_{i}^j - a_{i}^* \Vert <  \epsilon_0$ and $ \Vert a_{i-1}^* - a_{i-1}^j\Vert < \epsilon_0$.
			Therefore
			$\Vert  a_{i}^j - a_{i}^* \Vert < \Vert  c_{i} - a_{i}^*\Vert \text{ and } \Vert a_{i-1}^* - a_{i-1}^j\Vert < \Vert  c_{i-1} - a_{i-1}^*\Vert.$
			Besides,
			$ \epsilon_0 = \Vert  c_{i} - a_{i}^* \Vert =  \dfrac{\Vert s_2 - a_{i}^*  \Vert}{\Vert  s_2 - a_{i+1}^* \Vert} . \Vert  c_{i+1} - a_{i+1}^* \Vert = m_2. \Vert  c_{i+1} - a_{i+1}^* \Vert.$
			Combining with~(\ref{eq:parameter-M}) gives
			$	\Vert  a_{i+1}^j - a_{i+1}^* \Vert < \frac{\epsilon_0}{M}= \frac{m_2}{M}. \Vert  c_{i+1} - a_{i+1}^*\Vert \le \Vert  c_{i+1} - a_{i+1}^*\Vert.
			$
			These results yield the claim.
			% 			 Therefore
			% 			\begin{align}
			% 				\label{eq:claim-2a}
			% 				\Vert  a_{i}^j - a_{i}^* \Vert < \Vert  c_{i} - a_{i}^*\Vert \text{ and } \Vert a_{i-1}^* - a_{i-1}^j\Vert < \Vert  c_{i-1} - a_{i-1}^*\Vert.
			% 			\end{align}
			% 			Besides,
			% 			$$ \epsilon_0 = \Vert  c_{i} - a_{i}^* \Vert =  \dfrac{\Vert s_2 - a_{i}^*  \Vert}{\Vert  s_2 - a_{i+1}^* \Vert} . \Vert  c_{i+1} - a_{i+1}^* \Vert = m_2. \Vert  c_{i+1} - a_{i+1}^* \Vert.$$
			% 			Thus 
			% 			\begin{align}
			% 				\label{eq:claim-2b}
			% 				\Vert  a_{i+1}^j - a_{i+1}^* \Vert < \frac{\epsilon_0}{M}= \frac{m_2}{M}. \Vert  c_{i+1} - a_{i+1}^*\Vert \le \Vert  c_{i+1} - a_{i+1}^*\Vert.
			% 			\end{align} 
			% 			Combining~(\ref{eq:claim-2a}) with~(\ref{eq:claim-2b}) yields the claim.
		\end{proof}
		\begin{figure}[H]
			\centering
			\begin{subfigure}[b]{0.45\textwidth}
				%				\centering
				\raggedright
				\includegraphics[scale=0.52]{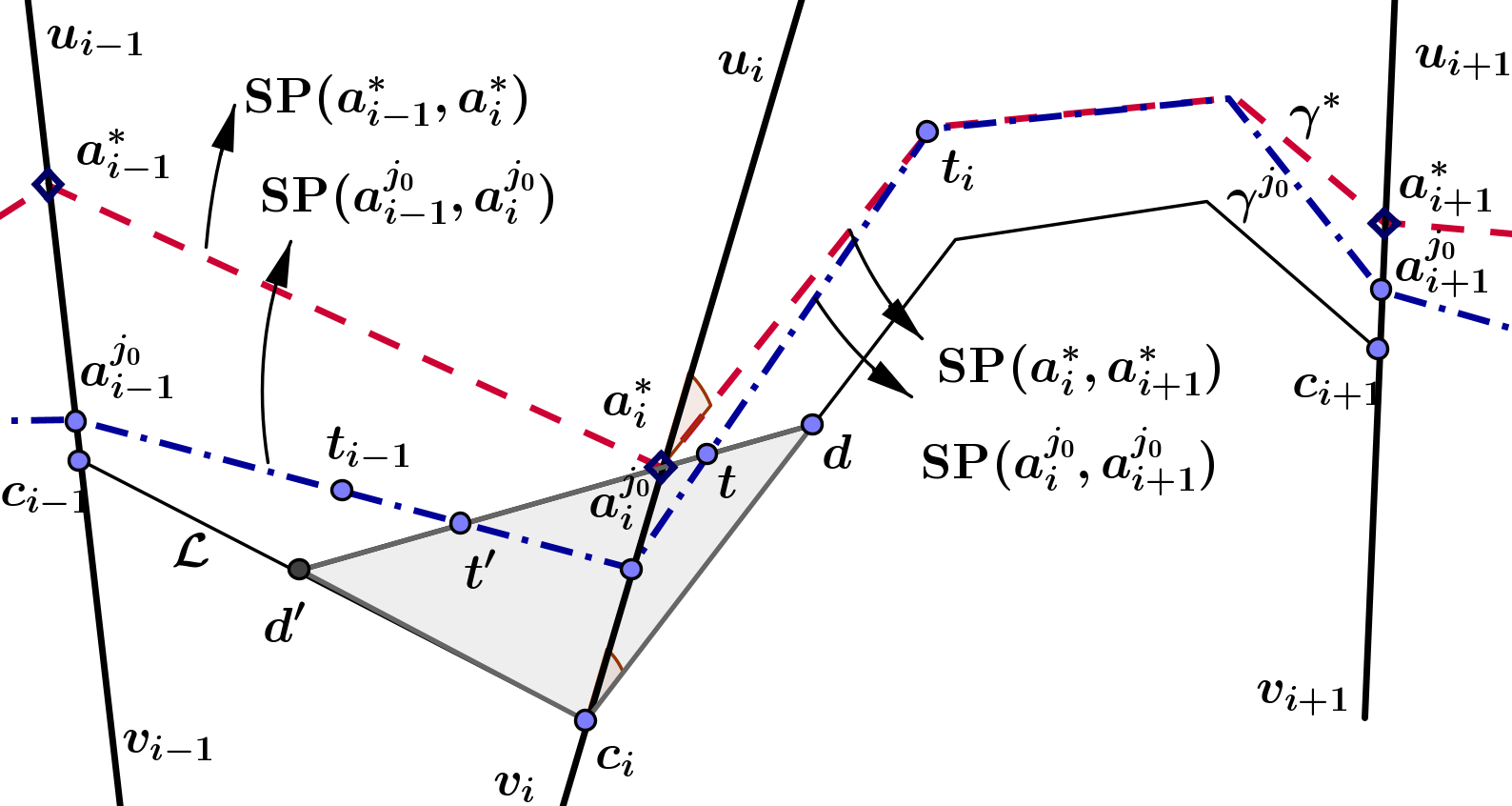}
				\caption*{(i)}
			\end{subfigure}
			\begin{subfigure}[b]{0.22\textwidth}
				\raggedleft
				\includegraphics[scale=0.46]{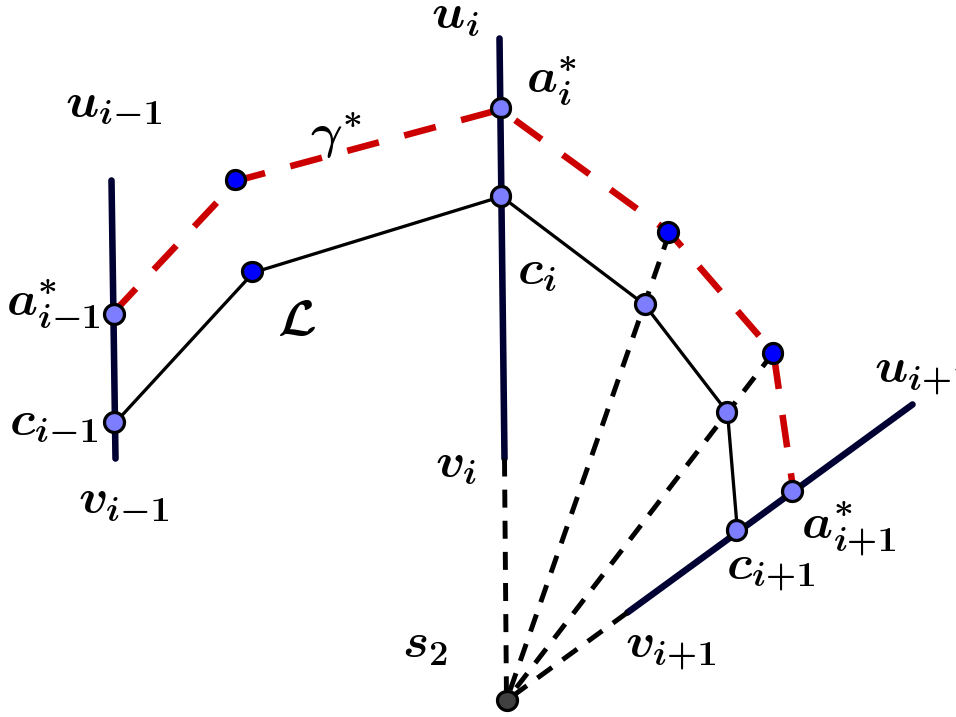}
				\caption*{(ii)}
			\end{subfigure}
			\begin{subfigure}[b]{0.31\textwidth}
				%	\centering
				\raggedleft
				\includegraphics[scale=0.6]{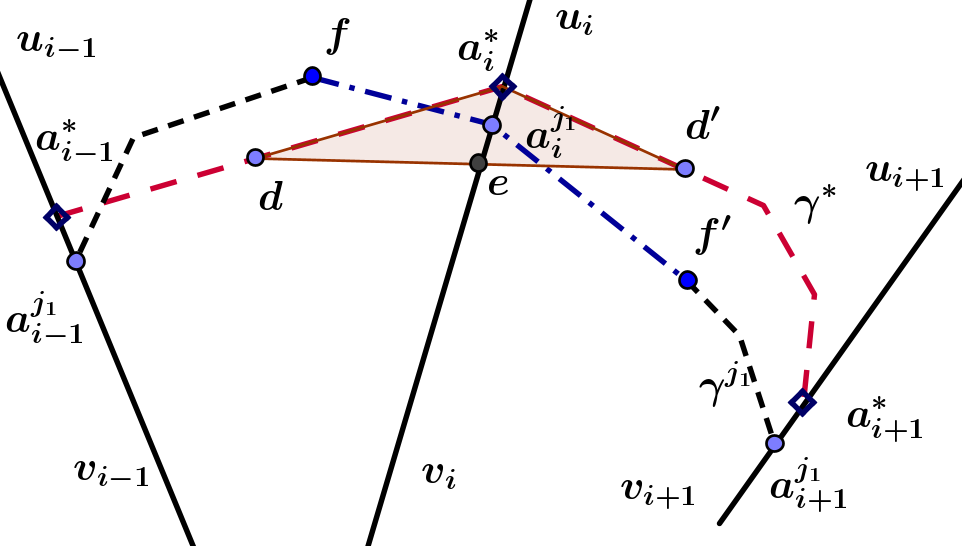}
				\caption*{(iii)}
			\end{subfigure}
			\caption{Illustration of the proof of Theorem~\ref{theo:global solution}}
			\label{fig:prooftheorem}
		\end{figure}
		\textit{We now turn to the proof of the theorem.}
		We have $\angle \left( \text{SP}(a_{i-1}^*,a_i^*), \text{SP}(a_i^*,a_{i+1}^*) \right) < \pi$ and the fact that $\{a_{i-1}^j\}$, $\{a_i^j\}$, and $\{a_{i+1}^j\}$  converge on one side to $a_{i-1}^*$, $a_i^*$ and $a_{i+1}^*$ as $j \rightarrow \infty$, respectively. There exists a natural number $j_0 \ge \overline{j}$  such that  $\angle \left( \text{SP}(a_{i-1}^{j_0},a_i^{j_0}), \text{SP}(a_i^{j_0},a_{i+1}^{j_0}) \right) < \pi$.
		Then the algorithm updates $a_i^{j_0}$ to $a_i^{j_0+1}$. \textit{We are now in a position to show $a_i^{j_0+1} \in (a_i^*, u_i]$.}

		Let two intersection points of $\gamma^{j_0}$ and $[d,d']$ be $t$ and $t'$.
		By Lemma~\ref{lem:AnHaiHoai}, \break $d_{\mathcal{H}}\left( {\rm SP}(a^{j_0}_i,a^{j_0}_{i+1}),{\rm SP} (a_i^*,a_{i+1}^*) \right)  \le 2 \max \{\Vert  a^{j_0}_{i} - a_{i}^*\Vert, \Vert  a^{j_0}_{i+1} - a_{i+1}^*\Vert \} 
		%\le \max \{\Vert  c_{i} - a_{i}^*\Vert, \Vert  c_{i+1} - a_{i+1}^*\Vert \}
		\le \mu/2.$
		Combining with~(\ref{eq:take-epsilon-0}), $\text{SP}(a_i^{j_0},a_{i+1}^{j_0})$ does not contain vertices of $V$ except the vertices belonging to $\hat{\gamma}$. Therefore, all endpoints  of  $\text{SP}(a_i^{j_0},a_{i+1}^{j_0})$ except for $a_i^{j_0}$ and $a_{i+1}^{j_0}$ belong to $\text{SP}(a_i^*,a_{i+1}^*)$.	
		Note that the temporary point $t_i$  of $\gamma^{j_0}$  is taken as in step~\ref{step:temp_points} of Procedure \uppercase{Collinear\_Update}($\gamma^{current}, \gamma^{next}, flag$). If $\text{SP}(a_i^{j_0},a_{i+1}^{j_0})$ is  a line segment, we have $\Vert a_i^{j_0} - t_i \Vert \ge \frac{1}{2} d([u_i,v_i],[u_{i+1},v_{i+1}]).$
		Otherwise, $t_i$ is one of the endpoints of $\text{SP}(a_i^{j_0},a_{i+1}^{j_0})$ except for $a_i^{j_0}$ and $a_{i+1}^{j_0}$, and thus $t_i \in \hat{\gamma}$. Hence  $\Vert a_i^{j_0} - t_i \Vert$ is greater than a given constant. 
		% We can assert that $\Vert a_i^{j_0} - t \Vert \le \Vert c_i - d \Vert$, and hence, by (\ref{eq:laws-of-sines}), it is very small positive numbers depending on $\epsilon_0$. 
		As $\epsilon_0$ is given by~(\ref{eq:arbitrarily small}), 		
		we have $t_i \notin [a_i^{j_0}, t]$. Similarly, we also get $t_{i-1}\notin [a_i^{j_0}, t']$. Because of the properties of shortest paths, 
		$\text{SP}(t_{i-1}, t_i)$ do not intersect with $[t,t']$. 		
		In the sequence,  $\text{SP}(t_{i-1}, t_i)$ intersects $[u_i,v_i]$ in a point that is above $a_i^*$ when we see from $u_i$. This contradicts the result proved above in which
		$\gamma^{j_0+1}$ is below $\hat{\gamma}$ when we view from $B_1$.

		\textbf{\textit{Case 2:}} $\angle \left( \text{SP}(a_{i-1}^*,a_i^*), \text{SP}(a_i^*,a_{i+1}^*) \right) > \pi$. Since $a_i^* \neq v_i$, there exist two points $d \in \text{SP}(a_{i-1}^*,a_i^*)$, and $d' \in \text{SP}(a_i^*,a_{i+1}^*)$ such that three points $d, d'$, and $a_i^*$  form a triangle that completely  belongs to   $\mathcal{D}$ but does not contain any vertex of $\mathcal{D}$, $[u_{i-1},v_{i-1}] \cap \triangle dd' a_i^* = \emptyset$, and $[u_{i+1},v_{i+1}]\cap \triangle dd' a_i^* = \emptyset$. Let $e = [d,d'] \cap [u_i, v_i]$.
		
		\textit{We claim that there is no element of $\{ a_i^j \}$ which is contained in $[e, a_i^*]$.}	
		Indeed, suppose that there exists  $j_1 \in \mathbb{N}$ such that $a_i^{j_1} \in [e,a_i^*]$. Let two segments having the common endpoint $a_i^{j_1}$ of $\gamma^{j_1}$ be $[f, a_i^{j_1}]$ and $[a_i^{j_1}, f']$, where $f \in  \text{SP}(a_{i-1}^{j_1},a_i^{j_1})$ and $f' \in \text{SP}(a_i^{j_1},a_{i+1}^{j_1})$.
		Then either $f=a_{i-1}^{j_1}$ or $f$ is  a vertex of $\mathcal{D}$, and hence $f \notin \triangle dd' a_i^*$.
		Similarly, we get $f' \notin  \triangle dd' a_i^*$.  According to Corollary~\ref{coro:inner-angle-shooting-point}, it follows that $\angle{fa_i^{j_1} f'} \le  \pi$. Because of $\angle d e d' = \pi$, we obtain either $[a_i^j, f]$ crossing $[a_i^*,d]$ or $[a_i^j,f']$ crossing $[a_i^*,d']$ (see Fig.~\ref{fig:prooftheorem}(iii)). This contradicts  the proved result that  $\gamma^{j_1}$ is below $\hat{\gamma}$ when we view from the boundary $B_1$. 
		%	\begin{figure}
		%		\centering
		%		\includegraphics[width=0.5\linewidth]{proof_theorem-5}
		%		\caption{Illustration of the proof of the theorem}
		%		\label{fig:prooftheorem-5}
		%	\end{figure}

		Summarizing, combining two cases gives that Collinear Condition (B) holds at $a_i^*$, for $i=1,2, \ldots, N$. Hence $\hat{\gamma}$ is the shortest path joining $\tilde{b}$ and $b$. The proof is complete.
	\end{proof}

	\medskip
	\noindent \textbf{The proof of Theorem~\ref{theo:update_path}:}
	
	\begin{proof}
		%	The proof can be handled in much the same way as that of Lemma 8.3.6 in~\cite{An2017}, the only difference 		being in the considered paths of the lemma which was in 3D.
		%\end{proof}
		
		%  Le 6/1/2022	
		(a) We have (see Figs.~\ref{fig:MMS}(iii) and (iv))
		\begin{align*}
			%	\label{eq:update_path_1}
			l(\gamma^{current}) = & \sum_{i=0}^{N}l(\text{SP}(a_i,a_{i+1})) \notag \\ 
			= & l(\text{SP}(a_0,t_0)) + \sum_{i=1}^{N} \left( l( \text{SP}(t_{i-1},a_i))+l(\text{SP}(a_i, t_{i})) \right) + l(\text{SP}(t_N, a_{N+1})) 
			\notag \\
			% 	\end{align*}
			%%
			%\begin{align*}
			\qquad	\qquad \qquad \quad
			\ge & \,  l(\text{SP}(a_0^{next},t_0)) + \sum_{i=1}^{N}  l \left( \text{SP}(t_{i-1},t_{i}) \right)  + l(\text{SP}(t_N, a_{N+1}^{next})) \notag \\
			= & \,  l(\text{SP}(a_0^{next},t_0)) + \sum_{i=1}^{N} \left( l(\text{SP}(t_{i-1},a_i^{next}))+l(\text{SP}(a_i^{next}, t_{i})) \right)  + l(\text{SP}(t_N, a_{N+1}^{next})) \notag \\
			= &  \sum_{i=0}^{N} \left( l(\text{SP}(a_i^{next},t_i)) +  l(\text{SP}(t_i,a_{i+1}^{next})) \right) \notag \\
			\ge & \sum_{i=0}^{N} \left( l (\text{SP}(a_i^{next},a_{i+1}^{next})) \right)
			=  l \left( \cup_{i=0}^{N} \text{SP}(a_i^{next},a_{i+1}^{next}) \right) = l (\gamma^{next}).
		\end{align*}
		
		If Collinear Condition (B) does not hold at some shooting point $a_i$, then first inequality is strict. Thus $l (\gamma^{current}) > l (\gamma^{next}).$
		
		(b)  Recall that $\gamma^* = \cup_{i=0}^N \text{SP}(a_i^*,a_{i+1}^*)$ is the shortest path joining $\tilde{b}$ and $b$ in $\mathcal{D}$ and
		$\gamma^j = \cup_{i=0}^N \text{SP}(a_i^j,a_{i+1}^j)$ is the sequence of paths obtained by the algorithm, where $a_i^*$ ($a_i^j$, resp.) is the intersection point between $\gamma^*$ ($\gamma^j$, resp.) and the corresponding cutting edge. If the algorithm stops after a finite number of iterations, the proof is trivial. Otherwise, repeat  the proof of Theorem~\ref{theo:global solution} (item \textbf{i}) we get $a^j_i \rightarrow a_i^*$ as $j \rightarrow \infty$, for $i=1,\ldots,N$, where $a^j_0 = a_0^*$ and $a^j_{N+1} = a_{N+1}^*$. According to~[\cite{HaiAn2011}, p. 542], 
		the geodesic distance between two points $x$ and $y$ in a simple polygon, which is measured by the length of the shortest path joining two points $x$ and $y$ in the simple polygon, is a metric on the simple polygon and it is continuous as a function
		of both $x$ and $y$. That is, $l({\rm SP}(a_i^j, a_{i+1}^j)) \rightarrow l({\rm SP}(a_i^*, a_{i+1}^*))$, for $i=0,1,\ldots,N$ as $j \rightarrow \infty$. It follows $l(\gamma^j) \rightarrow l(\gamma^*)$ as $j \rightarrow \infty$.
	\end{proof}
	
	\medskip

	\noindent \textbf{The proof of Theorem~\ref{theo:complexity}:}
	\begin{proof}
		%		Step 1 and 2 of the algorithm need $O(1)$ and $O(n)$ time, respectively. We know that the polygon $\mathcal{D}$ consists of $n+7$ vertices which are vertices of the polyline $\mathcal{P}$ and the rectangle $\mathcal{R}$. Besides, $\mathcal{D}$ is divided into $N+1$ sub-polygons $\mathcal{D}_i$.
		%		In step 3, taking  the set of initial shooting points takes $O(1)$ time.
		%Step 1-2 of the algorithm do not exceed $O(n)$ time.
		Steps~\ref{step1} and~\ref{step2} of the algorithm need $O(1)$ and $O(n)$ time, respectively.
		In step~\ref{step3}, there are at most two procedures of finding shortest paths in sub-polygons which call each vertex of $\mathcal{D}$. Since the problem of finding the shortest path joining two points in a simple polygon can be solved in linear time, computing the initial path $\gamma^{current}$ formed  by the set of initial shooting points takes $O(n)$ time.
		Steps~\ref{step6} and~\ref{step7} need $O(1)$ time. We need to show that step~\ref{step5} needs $O(n)$ time.
		%, i.e., the time complexity of Procedure \uppercase{Collinear\_Update}($\gamma^{current}, \gamma^{next}, flag$) is linear.
		%
		For each iteration step, 
		% since the geometry shortest path problem can be solved in linear time and there are at most two procedures of finding shortest paths in sub-polygons called each vertex of $\mathcal{D}$, finding the initial path $\gamma^{current}$ formed the set of initial shooting points takes $O(n)$ time.
		% each vertex of $\mathcal{D}$ is called at most twice in executing the procedure of finding shortest paths in sub-polygons). 
		there are $N$ cutting segments, where $N$ is a given constant number, %and $N+1$ sub-polygons $\mathcal{D}_i$, 
		then the computing sequentially the angles at shooting points is done in  constant time. In  step~\ref{step-prc8} and~\ref{step-prc14} of for-loop in the procedure,  since there are at most two procedures for finding the shortest path $\text{SP}(t_{i-1}, t_i)$ in sub-polygons $\mathcal{D}_{i-1} \cup \mathcal{D}_i$ which call each vertex of $\mathcal{D}$,
		% and the point of intersection of  $\text{SP}(t_{i-1}, t_i)$ and $[u_i, v_i]$ is found in 
		time complexity of Procedure \uppercase{Collinear\_Update}($\gamma^{current}, \gamma^{next}, flag$) is linear.
		In summarizing, the algorithm runs in $O(kn)$ time, where $k$ is the number of iterations to get the required path.
	\end{proof}


\begin{thebibliography}{50}
		
		\bibitem{ptan-convex-rope}  An, P. T. (2010). Reachable gasps on a polygon of a robot arm: 
		finding the convex rope without triangulation, \emph{Journal of Robotics and 
			Automation}, \textbf{25}(4),  304--310.
		
		\bibitem{An2017} An, P. T. (2017).  \textit{Optimization Approaches for Computational Geometry.} Publishing House for Science and Technology, Vietnam Academy of Science and Technology, Hanoi, ISBN 978-604-913-573-6.
		
		\bibitem{AnHaiHoai2013} An, P. T., Hai, N. N.,  Hoai, T. V.  (2013). Direct multiple shooting method
		for solving approximate shortest path problems. \textit{ Journal of Computational
			and Applied Mathematics}, \textbf{244},  67--76.
		
		\bibitem{AnHoai2010} An, P. T., Hoai, T. V. (2012). Incremental convex hull as an orientation to solving the shortest path problem. \textit{International Journal of Information and Electronics Engineering}, \textbf{2}(5), 652--655.
		
		\bibitem{AnTrang2018}  An, P. T., Trang., L. H. (2018).  Multiple shooting approach for computing shortest descending paths on convex terrains. \textit{Computational and Applied Mathematics}, \textbf{37}(4),  1--31. 
		
		\bibitem{AnHaiHoaiTrang2014}  An, P. T., Hai, N. N., Hoai, T. V., Trang., L. H. (2014).  On the performance of triangulation-based multiple shooting method for 2D geometric shortest path problems. \textit{Transactions on Large-Scale Data- and Knowledge-Centered Systems XVI},  45--56. 
		
		
		%
		
		%
		%\bibitem{bock-dmsm} Bock, H. G.,  Plitt, K. J. (1984). A multiple shooting algorithm for  direct solution optimal control problems, \emph{Proceedings of the 9th IFAC world  congress Budapest}, Pergamon Press,  243--247.
		%
		%	\bibitem{Aronov1993}  		Aronov, B., Fortune, S., Wilfong, G. (1993). The furthest-site geodesic Voronoi diagram. \textit{Discrete \& Computational Geometry}, \textbf{9}(3),  217--255.
		
		\bibitem{Chazelle1982} Chazelle, B. (1982). \textit{A theorem on polygon cutting with applications}, Proc. 23rd IEEE Symposium on Foundations of Computer Science, Chicago,  339--349.
		
		
		
		\bibitem{guibas-funnel} Guibas, L., Hershberger, J., Leven, D., Sharir, M., Tarjan, R. E. (1987). Linear-time algorithms for visibility and shortest path problems inside
		triangulated simple polygons, \emph{Algorithmica}, \textbf{2},  209--233.
		
		%
		\bibitem{HaiAn2011} Hai, N. N.,  An, P. T. (2011). Blaschke-type theorem and separation of disjoint closed geodesic convex sets. \textit{Journal of Optimization Theory and Applications}, \textbf{151}(3), 541-551.
		
			\bibitem{HaiAnHuyen2019} Hai, N. N., An, P. T., \& Huyen, P. T. T. (2019). \textit{Shortest paths along a sequence of line segments in euclidean spaces}. Journal of Convex Analysis, \textbf{26}(4), 1089-1112.
			
		\bibitem{heffernan} Heffernan, P. J.,  Mitchell, J. S. B. (1990). Structured visibility profiles 		with applications to problems in simple polygons, \emph{Proc. 6th Annual ACM Symp. 			Computational Geometry},   53--62.
		
		\bibitem{HoaiAnHai2017} Hoai, T. V. , An, P. T., Hai, N. N. (2017).  Multiple shooting approach for computing approximately shortest paths on convex polytopes. \textit{Journal of Computational and Applied Mathematics}, \textbf{317},  235--246.
		
		
		
		\bibitem{Lee1984} Lee, D. T., Preparata, F. P. (1984). Euclidean shortest paths in the presence of rectilinear
		barriers, \textit{Networks}, \textbf{14}(3),  393--410.
		
		
		\bibitem{Li2006}  Li. F,  Klette. R. (2006). Finding the shortest path between two points in a simple polygon by applying a rubberband algorithm. \textit{In Pacific-Rim Symposium on Image and Video Technology}. Springer,  Heidelberg,  280--291.
		%
		\bibitem{perez} Lozano-Perez, T.  (1981). Automatic planning of manipulator transfer movements,
		\emph{IEEE Transactions on Systems, Man, and Cybernetics}, SMC-11(10), 
		681-698.
		%
		\bibitem{melkman}  Melkman, A. A. (1987). On-line construction of the convex hull of a simple
		polyline, \emph{Information Processing Letters}, \textbf{25},  11--12.
		%
		%\bibitem{mitchell} J. S.B. Mitchell, Geometric shortest paths and network optimization,
		%in \emph{Handbook of Computational Geometry}, J.. R. Sack and J. Urrutia, eds, Elsevier
		%Science B. V., 2000,  633-701.
		%
		%\bibitem{orourke} J. O'Rourke, Computational Geometry in C, Cambridge University Press, 
		%Second Edition, 1998.
		%
		%	\bibitem{Papadopoulos2005} Papadopoulos, A. (2005). \textit{Metric spaces, convexity and nonpositive curvature}, European Mathematical Society.
		
		\bibitem{peshkin} Peshkin, M. A., Sanderson, A. C. (1986). Reachable grasps on a polygon: 
		the convex rope algorithm, \emph{IEEE Journal on Robotics and Automation}, \textbf{2}(1),  5--58.
		%
		
		
		\bibitem{SongAtAl2018} Song, P., Fu, Z., and Liu, L. (2018). Grasp planning via hand-object geometric fitting. \textit{The Visual Computer}, \textbf{34}(2), 257--270.
		
		
		\bibitem{Toussaint1989} Toussaint, G. T. (1989). Computing geodesic properties inside a simple polygon, \textit{Revue d'Intelligence Artificielle}, \textbf{3}(2),  9--42.
		
		
		
		
		\bibitem{Xue2008} Xue, Z., Zoellner, J. M.,  Dillmann, R.  (2008).  \textit{Automatic optimal grasp planning based on found contact points}. International Conference on Advanced Intelligent Mechatronics, IEEE,  1053--1058.
		
		%
	\end{thebibliography}
\end{document}